\documentclass[12pt]{amsart}
\usepackage{amsmath,amscd,amssymb,amsfonts,graphics}
\usepackage{xcolor,hyperref}
\hypersetup{colorlinks=true,citecolor=black,linkcolor=black}
\newcommand{\hl}{\hyperlink}
\newcommand{\htt}{\hypertarget}
\newcommand{\h}{\hbox}
\newcommand{\q}{\quad}
\newcommand{\nin}{\noindent}
\newcommand{\bs}{\par\bigskip}
\newcommand{\ms}{\par\medskip}
\newcommand{\sk}{\par\smallskip}
\newcommand{\bsn}{\par\bigskip\noindent}
\newcommand{\msn}{\par\medskip\noindent}
\newcommand{\skn}{\par\smallskip\noindent}
\newcommand{\ssb}{\raise.15ex\h{${\scriptscriptstyle\bullet}$}}
\newcommand{\ssc}{\,\raise.15ex\h{${\scriptstyle\circ}$}\,}
\newcommand{\msum}{\hbox{$\sum$}}
\newcommand{\mprod}{\hbox{$\prod$}}
\newcommand{\mcup}{\hbox{$\bigcup$}}
\newcommand{\mcap}{\hbox{$\bigcap$}}

\newcommand{\mopl}{\hbox{$\bigoplus$}}
\newcommand{\C}{{\mathbb C}}
\newcommand{\DD}{{\mathbf D}}
\newcommand{\N}{{\mathbb N}}
\newcommand{\PP}{{\mathbb P}}
\newcommand{\Q}{{\mathbb Q}}
\newcommand{\R}{{\mathbf R}}
\newcommand{\Z}{{\mathbb Z}}
\newcommand{\A}{{\mathcal A}}
\newcommand{\B}{{\mathcal B}}
\newcommand{\Cc}{{\mathcal C}}
\newcommand{\D}{{\mathcal D}}
\newcommand{\F}{{\mathcal F}}
\newcommand{\G}{{\mathcal G}}
\newcommand{\Hc}{{\mathcal H}}
\newcommand{\I}{{\mathcal I}}
\newcommand{\K}{{\mathcal K}}
\newcommand{\Lc}{{\mathcal L}}
\newcommand{\M}{{\mathcal M}}
\newcommand{\Nc}{{\mathcal N}}
\newcommand{\OO}{{\mathcal O}}
\newcommand{\naal}{\nabla^{(\alpha)}}
\newcommand{\alt}{\widetilde{\alpha}}
\newcommand{\jt}{{}\,\widetilde{\!j}{}}
\newcommand{\bt}{\widetilde{b}}

\newcommand{\Lt}{\widetilde{L}}
\newcommand{\Dt}{\widetilde{D}}
\newcommand{\Vt}{\widetilde{V}}
\newcommand{\Xt}{\widetilde{X}}

\newcommand{\pit}{\widetilde{\pi}}
\newcommand{\xit}{\widetilde{\xi}}
\newcommand{\Tht}{\widetilde{\Theta}}
\newcommand{\mm}{{\mathfrak m}}
\newcommand{\al}{\alpha}
\newcommand{\la}{\lambda}

\newcommand{\Ga}{\Gamma}
\newcommand{\dd}{\partial}
\newcommand{\ddd}{{\rm d}}
\newcommand{\om}{\omega}
\newcommand{\Om}{\Omega}
\newcommand{\DR}{{\rm DR}}
\newcommand{\DRl}{{}^l\hskip-.7pt{\rm DR}}
\newcommand{\Gr}{{\rm Gr}}
\newcommand{\bl}{\bigl}
\newcommand{\br}{\bigr}
\newcommand{\into}{\hookrightarrow}
\newcommand{\simto}{\,\rlap{\hskip1.6mm\raise1.5mm\hbox{$\sim$}}\hbox{$\longrightarrow$}\,}
\newcommand{\simfrom}{\,\rlap{\hskip1.9mm\raise1.5mm\hbox{$\sim$}}\hbox{$\longleftarrow$}\,}
\newcommand{\onto}{\twoheadrightarrow}
\newcommand{\ontov}{\rotatebox{-90}{$\onto$}}
\newcommand{\ges}{\geqslant}
\newcommand{\les}{\leqslant}
\newcommand{\1}{\hskip1pt}
\newcommand{\sst}{\,{\subset}\,}
\newcommand{\stm}{\,{\setminus}\,}
\newcommand{\pl}{\1{+}\1}
\newcommand{\mi}{\1{-}\1}
\newcommand{\eq}{\,{=}\,}
\newcommand{\nes}{\,{\ne}\,}
\newcommand{\less}{\,{\les}\,}
\newcommand{\gess}{\,{\ges}\,}
\newcommand{\ins}{\,{\in}\,}
\newcommand{\tos}{\,{\to}\,}
\newcommand{\caps}{\1{\cap}\1}
\newcommand{\indlim}{\rlap{\raise-6pt\h{$\,\to$}}{\rm lim}}
\begin{document}
\title[Twisted logarithmic complexes]{Twisted logarithmic complexes of positively\\weighted homogeneous divisors}
\author{Daniel Bath}
\address{Department of Mathematics, KU Leuven, Celestijnenlaan 200B, 3001 Leuven, Belgium}
\email{dan.bath@kuleuven.be}
\author{Morihiko Saito}
\address{RIMS Kyoto University, Kyoto 606-8502 Japan}
\email{msaito@kurims.kyoto-u.ac.jp}
\thanks{The first named author is supported by FWO grant \#G097819N and FWO grant \#12E9623N. The second named author was partially supported by JSPS Kakenhi 15K04816.}
\begin{abstract} For a rank 1 local system on the complement of a reduced divisor on a complex manifold $X$, its cohomology is calculated by the twisted meromorphic de Rham complex. Assuming the divisor is everywhere positively weighted homogeneous, we study necessary or sufficient conditions for a quasi-isomorphism from its twisted logarithmic subcomplex, called the logarithmic comparison theorem (LCT), by using a stronger version in terms of the associated complex of $D_X$-modules. In case the connection is a pullback by a defining function $f$ of the divisor and the residue is $\alpha$, we prove among others that if LCT holds, the annihilator of $f^{\alpha-1}$ in $D_X$ is generated by first order differential operators and $\alpha{-}1{-}j$ is not a root of the Bernstein-Sato polynomial for any positive integer $j$. The converse holds assuming either of the two conditions in case the associated complex of $D_X$-modules is acyclic except for the top degree. In the case where the local system is constant, the divisor is defined by a homogeneous polynomial, and the associated projective hypersurface has only weighted homogeneous isolated singularities, we show that LCT is equivalent to that $-1$ is the unique integral root of the Bernstein-Sato polynomial. We also give a simple proof of LCT in the hyperplane arrangement case under appropriate assumptions on residues, which is an immediate corollary of higher cohomology vanishing associated with Castelnuovo-Mumford regularity. Here the zero-extension case is also treated.
\end{abstract}
\maketitle
\centerline{\bf Introduction}
\bsn
Let $D$ be a reduced divisor on a complex manifold $X$ of dimension $n$. Let $L$ be a rank 1 local system on the complement $U:=X\stm D$ with $\la_k$ eigenvalues of the local monodromies of $L$ around each global irreducible component $D_k$ of $D$. Choosing complex numbers $\al_k$ with $e^{-2\pi i\1\al_k}\eq\la_k$, we have a locally free $\OO_X$-module $\Lc_X$ of rank 1 endowed with a meromorphic integrable connection $\naal$ which has a pole along $D$ and calculates the local system $L$ on $U$. Choosing local defining functions $f_k$ of $D_k\sst X$, it can be defined {\it locally\1} by using the {\it twisted differential}
\htt{1}{}
$$\ddd\pl\msum_k\,\al_k\1\om_k\wedge:\OO_X\to\Om_X^1(\log D)\q\h{with}\q\om_k:=\ddd f_k/\!f_k,
\leqno(1)$$
\par\nin and trivializing $\Lc_X$, where $\Om_X^j(\log D)$ is the sheaf of logarithmic forms, see \cite{SaK}. We can apply the Hartogs extension theorem to see the independence of $\Lc_X$. (This is well known at the normal crossing locus of $D$, see \cite{De}.)
\sk
Setting $\M_X(L):=\Lc_X(*D)$, we get the {\it meromorphic de Rham complex}
$$(\M_X^{\ssb}(L),\naal)\,\bl(=\DR_X(\M_X(L),\naal)[-n]\br),$$
\par\nin which is the de Rham complex associated with a meromorphic connection in the classical sense, see \cite{De}. This is isomorphic to $\R(j_U)_*L$ in the derived category $D^b_c(X,\C)$ with $j_U:U\into X$ the inclusion (by reducing to the normal crossing case, see \cite{De}). Note that $\M_X^0(L)=\M_X(L)$, since the de Rham functor $\DR_X$ is shifted by $n$ as usual.
\sk
Employing the logarithmic differential forms, we can get the {\it logarithmic subcomplex}
\htt{2}{}
$$\bl(\M^{\ssb}_{X,\log}(L),\naal\br)\into\bl(\M_X^{\ssb}(L),\naal\br),
\leqno(2)$$
\par\nin with $\M^j_{X,\log}(L):=\Om_X^j(\log D){\otimes}_{\OO_X}\Lc_X$ for $j\ins\Z$. It is interesting whether (\hl{2}{2}) is a {\it quasi-isomorphism.} There is a lot of work about this problem, see for instance \cite{CNM}, \cite{WY}, \cite{HM}, \cite{To2}, \cite{CN}, \cite{Ba} (and also \S\hl{S3}{3} below). In the {\it plane curve\1} case with $L\eq\C_U$, the quasi-isomorphism is equivalent to {\it positively weighted local homogeneity\1} of $D$, see \cite{CMNC}. We can extend the assertion to the higher dimensional case assuming {\it isolated singularities,} see Corollary~\hl{C1}{1} below. This fails without assuming the latter condition, for instance, in the case $f\eq xy(x{+}y)(x{+}yz)$ where the logarithmic comparison theorem holds rather surprisingly, see \cite[\S3.1]{Na1} (and \cite[\S4]{CMNC} on $\{z(z{-}1)\nes0\}\sst\C^3$).
\sk
Since (\hl{2}{2}) is a morphism of {\it differential complexes\1} in the sense of \cite{mhp}, \cite{ind}, we may consider a stronger version: Is (\hl{2}{2}) a $D$-quasi-isomorphism? Here a $D$-quasi-isomorphism means that we get a quasi-isomorphism applying the functor $\DR_X^{-1}$ to the morphism (\hl{2}{2}), see \cite[(1.3.2)]{ind}. (This is defined by $\otimes_{\OO_X}\D_X$ for each component of the logarithmic complex, and has nothing to do with Kashiwara's construction of the inverse functor of the Riemann-Hilbert correspondence, where only the underlying constructible $\C$-complexes are considered forgetting the $\OO_X$-module structure.)
This stronger version is equivalent to the canonical quasi-isomorphism
\htt{3}{}
$$\DRl_X^{-1}\bl(\M^{\ssb}_{X,\log}(L),\naal\br)[n]\simto\M_X(L).
\leqno(3)$$
\par\nin Here $\DRl_X^{-1}$ is the composition of $\DR_X^{-1}$ with ${\otimes}_{\OO_X}\om_X^{\vee}$, and the latter is the inverse of the transformation from left $\D_X$-modules to the corresponding right $\D_X$-modules. (For the advantage of considering this stronger version, see for instance Theorem~\hl{T2}{2} below.) We can verify that the above two questions are {\it equivalent\1} to each other in the everywhere positively weighted homogeneous (or more generally, locally finite logarithmic stratification) case, see Lemma~\hl{L1.6}{1.6} and Proposition~\hl{P1.6}{1.6} below. We can solve these equivalent questions in the hyperplane arrangement case as an immediate corollary of higher cohomology vanishing associated with the Castelnuovo-Mumford regularity of logarithmic forms, which can be shown in the same way as in the logarithmic vector field case (\cite{Sc}, \cite{DeSi}, \cite{4var}), see Corollaries~\hl{C3.1}{3.1} and \hl{C3.2a}{3.2a}-\hl{C3.2b}{b} below (and also Corollary~\hl{C3.2c}{3.2c} for the {\it zero-extension\1} case).
\sk
We say that a reduced divisor $D$ is {\it positively weighted homogeneous\1} around $p\ins D$ if there are local coordinates $x_1,\dots,x_n$ of $X$ with center $p$ such that $D\sst X$ is locally defined by a weighted homogeneous polynomial $f_p$ of strictly positive weights $w_{p,i}$, that is, $f_p$ is a linear combination of monomials $\mprod_{i=1}^n\,x_i^{\nu_i}$ with $\msum_{i=1}^n\,w_{p,i}\nu_i\eq 1$. We say that $D$ is {\it everywhere positively weighted homogeneous\1} if the above condition is satisfied at any $p\in D$. (This is called locally (or strongly) quasi-homogeneous in \cite{CNM}, \cite{CN}, \cite{Na2}.)
\sk
There is an irreducible factorization $f_p\eq\mprod_k\,f_{p,k}$, and we have $\al_{p,k}\ins\C$ by choosing globally the $\al_k$ as above. Put
\vskip-7mm
$$d_{p,k}:=\deg_{w_p}f_{p,k}.$$
\par\nin Here $\deg_{w_p}$ is the weighted degree associated with the weights $w_{p,i}$ of variables $x_i$, and is defined by the eigenvalue of the weighted Euler field $\msum_i\,w_{p,i}x_i\dd_{x_i}$. So $\msum_k\,d_{p,k}\eq 1$. (In the hyperplane arrangement case, $d_{p,k}\eq(\deg f_p)^{-1}$.) Set
\htt{4}{}
$$\alt_p:=\msum_k\,d_{p,k}\1\al_{p,k},\q e_p:=\min\bl\{e\ins\Z_{>0}\mid e\1w_{p,i}\ins\Z\br\}.
\leqno(4)$$
\par\nin These are locally constant on each stratum of the logarithmic stratification (which is locally finite in this case), see \hl{1.6}{1.6} below. In this paper we prove the following.
\par\htt{T1}{}\msn
{\bf Theorem~1.} {\it Assume $D$ is positively weighted homogeneous around $p$. Then the stalks at $p$ of the source and target of $(\hl{2}{2})$ are both acyclic $($hence $(\hl{2}{2})$ is a quasi-isomorphism at $p\1)$ if}
\htt{5}{}
$$e_p\1\alt_p\notin\Z.
\leqno(5)$$
\par\nin \sk
This follows from a calculation of the twisted logarithmic complex in Proposition~\hl{P1.5}{1.5}, see \hl{2.1}{2.1} below. Some more information is available in special cases of free divisors or hyperplane arrangements, see \cite{CNM}, \cite{CN}, \cite{Ba} (and also \S\hl{S3}{3} below).
\sk
We then prove the following, which is very much inspired by Theorem~\hl{TA}{A} in Appendix, and gives a partial answer to a generalization of a question in \cite[3.3]{To2}.
\par\htt{T2}{}\msn\vbox{\nin
{\bf Theorem~2.} {\it Assume $\al_k\eq\al$ $(\forall\,k)$ for some $\al\ins\C$. Consider the following conditions\,$:$
\par\htt{a2}{}\skn
{\rm(a)} The comparison morphism~{\rm(\hl{2}{2})} is a $D$-quasi-isomorphism.
\par\htt{b2}{}\skn
{\rm(b)} The comparison morphism~{\rm(\hl{2}{2})} is a quasi-isomorphism.
\par\htt{c2}{}\skn
{\rm(c)} The annihilator ${\rm Ann}_{\D_X}(f^{\al-1})$ is generated by $\Tht_{f,\al-1}$, see Remark~{\rm\hl{R1.7a}{1.7a}} below.
\par\htt{d2}{}\skn
{\rm(d)} The annihilator ${\rm Ann}_{\D_X}(f^{\al-1})$ is generated by first order differential operators.
\par\htt{e2}{}\skn
{\rm(e)} We have $b_f(\al{-}1{-}j)\nes0$ for any $j\ins\Z_{>0}$.
\par\htt{f2}{}\skn
{\rm(f)} The $\D_X$-module $\OO_X(*D)f^{\al-1}$ is generated by $f^{\al-1}$.
\par\skn
In general, condition~{\rm(a)} implies the other conditions, and we have the equivalences {\rm(c)}\,$\Leftrightarrow$\,{\rm(d)} and {\rm(e)}\,$\Leftrightarrow$\,{\rm(f)} unconditionally. 
If $D$ is everywhere positively weighted homogeneous and is tame $($or more generally, $\DRl_X^{-1}\bl(\M^{\ssb}_{X,\log}(L),\naal\br)[n]$ is quasi-isomorphic to a $\D_X$-module, see Corollary~{\rm\hl{C1.7b}{1.7b}} below$)$, then the above six conditions are equivalent to each other.}}
\ms
Here we assume the existence of Bernstein-Sato polynomial $b_f(s)$, which is called the BS polynomial for short in this paper, shrinking $X$ if necessary. For the definition of tame, see Corollary~{\rm\hl{C1.4}{1.4}} below. The equivalence {\rm(e)}\,$\Leftrightarrow$\,{\rm(f)} follows for instance from \cite[Thm.\,1]{rp}, see Proposition~\hl{P2.2}{2.2} below. In case the two assumptions on $D$ for the last assertion are satisfied, we can show {\rm(c)}\,$\Rightarrow$\,{\rm(a)} and {\rm(f)}\,$\Rightarrow$\,{\rm(a)} as in \hl{2.2}{2.2} below. We also give a proof of {\rm(c)}\,$\Leftrightarrow$\,{\rm(e)} using a quite different method in \hl{Ap}{Appendix} (where the tameness assumption cannot be weakened). It is not necessarily easy to construct an example where one can prove that $\DRl_X^{-1}\bl(\M^{\ssb}_{X,\log}(L),\naal\br)[n]$ is {\it not\1} quasi-isomorphic to a $\D_X$-module (even in the non-free divisor case). This is partly related to the {\it torsion-freeness\1} of the lowest degree (possibly) non-zero part of the Brieskorn modules (see for instance \cite[(4.3.6)]{DS2}). Such an example is given by $f\eq(xz{+}\1y)(x^4{+}\1y^5{+}\1xy^4)$, see \cite[Ex.\,5.1]{CN0} (and Remark~\hl{R1.7e}{1.7e} below).
\sk
From now on, we assume in the introduction the following:
\htt{6}{}
$$\h{$L\eq\C_U$, and the residues of connection $\al_k$ are 0.}
\leqno(6)$$
\par\nin In the {\it isolated\1} singularity case, combining Theorem~2 for $\al\eq0$ with \cite[Thm.\,1.2]{To1} (and also Proposition~\hl{P1.6}{1.6} below), we get the following.
\par\htt{C1}{}\msn
{\bf Corollary~1.} {\it Let $D\sst X$ be a hypersurface having only isolated singularities. Assume the logarithmic comparison theorem holds, that is, the comparison morphism~{\rm(\hl{2}{2})} for $\al\eq0$ is a quasi-isomorphism. Then $D$ is everywhere positively weighted homogeneous.}
\ms
Its converse holds if $n\eq 2$, and assuming the conclusion of Corollary~\hl{C1}{1} in the case $n\gess 3$, the logarithmic comparison theorem holds if and only if $D$ is a $\Q$-homology manifold, see \cite{HM} (and also Proposition~\hl{P1.4}{1.4} below). Here we use the Wang sequence \cite {Mi} (as is explained below) together with the {\it symmetry of spectral numbers,} see \cite{St1}, \cite{St2} (and also \cite{JKSY}). Corollary~\hl{C1}{1} has been shown in \cite[Thm.\,2]{Scl} under certain hypotheses, for instance, in the case where the minimal spectral number is greater than 1, that is, $D$ has only {\it rational singularities\1} \cite{exp}, or 1 is not an eigenvalue of the monodromy, that is, $D$ is a $\Q$-{\it homology manifold\1} (using the Wang sequence \cite{Mi}). The last conditions of (c) and (d) in \cite[Thm.\,2]{Scl} are not always satisfied, and do not seem easy to verify in general, see for instance \cite{sem}. Note that Corollary~\hl{C1}{1} implies a similar assertion around {\it general\1} points of ${\rm Sing}\,D$ if the logarithmic stratification is {\it locally finite,} for instance, if $D$ is defined by a positively weighted homogeneous polynomial and $\dim{\rm Sing}\,D\eq1$, see Remark~\hl{R1.6a}{1.6a} below.
\sk
Returning to the {\it general\1} singularity case, consider the morphisms
\htt{7}{}
$$H^j(\Om_{X,p}^{\ssb}(\log D),\ddd)\to H^j(\Om_{X,p}^{\ssb}(*D),\ddd)\q(j\ins\N,\,p\ins D),
\leqno(7)$$
\par\nin where $(*D)$ denotes the localization along $D$. (In the hyperplane arrangement case, these are surjective as a consequence of \cite[Lem.\,5]{Br2}.) We can handle the {\it weighted homogeneous\1} case much better because of the key Proposition~\hl{P1.4}{1.4} below, which shows {\it a close relation between logarithmic complexes and Brieskorn modules\1} \cite{Br}, \cite{BS} in the weighted homogeneous case. We can prove for instance the following.
\par\htt{T3}{}\msn\vbox{\nin
{\bf Theorem~3.} {\it Assume $D\sst X$ is defined by a positively weighted homogeneous polynomial $f$ around $p\ins D$.
\skn
{\rm(i)} For $j\ins\Z$, the morphism {\rm(\hl{7}{7})} is injective if and only if the weighted degree $1$ part of Brieskorn modules $\bl(\Hc_{f,p}^j\br)_1$ and $\bl(\Hc_{f,p}^{j+1}\br)_1$ has no torsion, see {\rm\hl{1.2}{1.2}} below for the notation.
\skn
{\rm(ii)} Assume that $-1$ is the unique integral root of the local BS polynomial $b_{f,q}(s)$ at $q\ne p$ sufficiently near $p$. Then the morphism {\rm(\hl{7}{7})} for $j\eq n$ is surjective if and only if $-1$ is the unique integral root of the BS polynomial $b_f(s)$.}}
\ms
Recall that the roots of $b_{f,q}(s)$ and $b_f(s)$ are {\it strictly negative\1} rational numbers, see \cite{Ka}. Theorem~\hl{T3}{3} follows from the coincidence of the torsion of Brieskorn modules with the kernel of the morphism to the Gauss-Manin systems (see \cite[Thm.\,1]{BS} or \hl{1.2}{1.2} below) together with the relation between the pole order filtrations \cite{DS2} and the roots of BS polynomials (see \cite[Thm.\,2]{bCM} or Theorem~\hl{T1.3}{1.3} below), where the key Proposition~\hl{P1.4}{1.4} below is used in an essential way, see \hl{2.3}{2.3} below for details.
\sk
Note that the Brieskorn modules $\Hc_{f,p}^j$ are {\it torsion-free\1} under the actions of $t$ and $\dd_t^{-1}$ for $1\,{<}\,j\less{\rm codim}_X\1{\rm Sing}\,f$, and vanish for $1\,{<}\,j\,{<}\,{\rm codim}_X\1{\rm Sing}\,f$ by the same argument as in \cite[(4.3.6)]{DS2}. Its corresponding Milnor fiber cohomology vanishes as a consequence of the semi-perversity of vanishing cycle complexes, see \cite{DS1}. Here we can use also the microlocal Gauss-Manin system as in \cite{DS2} together with \cite[Prop.\,2]{JKSY2}. For $j\eq 1$, $\Hc_{f,p}^1$ is a {\it free\1} $\C\{t\}$-module of rank 1, which is generated by $\ddd f$, and is isomorphic to $\C\{t\}$ endowed with the natural action of $\dd_t$, that is, $\dd_t[\ddd f]\eq0$ (since $f$ is reduced), see for instance \cite[\S 2]{BS}. We have the following.
\par\htt{C2}{}\msn
{\bf Corollary~2.} {\it Assume that $X\eq\C^n$ with $n\gess 3$, $D$ is defined by a homogeneous polynomial $f$, the projective hypersurface $Z\sst\PP^{n-1}$ defined by $f$ has only isolated singularities, and moreover these are all weighted homogeneous. Then the morphism~{\rm(\hl{7}{7})} is injective for any $j\ins\Z$ and $p\ins D$. It is surjective for any $j\ins\Z$ and $p\ins D$ if and only if $-1$ is the unique integral root of $b_f(s)$.}
\ms
This is a corollary of Theorem~\hl{T3}{3} combined with \cite[Thm.\,2]{wh} showing that the {\it pole order spectral sequence\1} degenerates at $E_2$  (or equivalently, the Brieskorn module $\Hc_{f,0}^n$ is {\it torsion-free,} see \cite[Cor.\,4.7]{DS2}) under the hypothesis of Corollary~\hl{C2}{2}. Since $b_f(s)$ is divisible by $b_{f,p}(s)$ for $p\nes0$, the condition about roots of $b_f(s)$ in Corollary~\hl{C2}{2} implies that (\hl{2}{2}) induces an isomorphism on $X\stm\{0\}$. For $\Hc_{f,0}^{n-1}$, we can use \cite[Thm.\,0.1]{DS1} (see Remark~\hl{R1.3c}{1.3c} below) together with Remark~\hl{R1.4a}{1.4a} below and also the inclusion $F\sst P$ (see \cite[Prop.\,4.4]{bCM}) for the case $n\eq 3$. Under the hypotheses of Corollary~\hl{C2}{2}, we can show that $\DRl_X^{-1}\bl(\Om_X^{\ssb}(\log D),\ddd\br)[n]$ is quasi-isomorphic to a {\it regular holonomic $\D_X$-module\1} (which is not necessarily a submodule of $\OO_X(*D)$), see Remark~\hl{R1.7c}{1.7c} below. Note also that $-n/d$ is not necessarily a roof of $b_f(s)$ under the hypotheses of Corollary\hl{C2}{2}, for instance if $f\eq x_1^ax_2^b\pl\msum_{i=3}^n\,x_i^d$ for $a,b\gess 2$ with $a,b$ mutually prime and $a{+}b\eq d$.
\sk
For $n\gess 4$, the morphism (\hl{7}{7}) is {\it not\1} surjective under the hypothesis of Corollary~\hl{C2}{2} in almost all cases. There is, however, a {\it quite exceptional case\1} for $n\eq 4$, where the morphism~(\hl{2}{2}) is a quasi-isomorphism with assumptions of Corollary~\hl{C2}{2} all satisfied; in particular, $-1$ is the unique integral root of $b_f(s)$. This happens if $f\eq x^d\pl g(y,z)\1w$ with $g(y,z)$ a reduced homogeneous polynomial of degree $d{-}1$ in $y,z$ and $d\gess 3$ (for instance, $g\eq y^{d-1}\pl z^{d-1}$), where \cite[Rem.\,3.2]{low} and the Thom-Sebastiani type theorem \cite[Thm.\,0.8]{mic} are used. It is unknown if there are other such examples for $n\gess 4$. If $n\eq 3$, there are many (see for instance \cite{DiSt1}, \cite{wh}), and we have the following.
\par\htt{P1}{}\msn\vbox{\nin
{\bf Proposition~1.} {\it Under the notation and assumptions of Corollary~{\rm\hl{C2}{2}}, assume further $n\eq3$. Then the following conditions are equivalent\,$:$
\par\htt{a}{}\skn
{\rm(a)} The morphism~{\rm(\hl{2}{2})} is a quasi-isomorphism on $X$.
\par\htt{b}{}\skn
{\rm(b)} The morphism~{\rm(\hl{7}{7})} for $j\eq 3$ is surjective at $0\ins X$.
\par\htt{c}{}\skn
{\rm(c)} $\,-1$ is the unique integral root of $b_f(s)$.
\par\htt{d}{}\skn
{\rm(d)} $H^0_{\mm}\bl(\OO_{X,0}/(\dd f)\br){}_{d-3}\eq 0$, or equivalently, $M'_d\eq0$, see {\rm(\hl{2.4.3}{2.4.3})} below.
\par\htt{e}{}\skn
{\rm(e)} All the roots of $b_f(s)$ are contained in $(-2,0)$.
\par\skn
Here $\mm$ and $(\dd f)$ are respectively the maximal and Jacobian ideals of $\OO_{X,0}$, and $d\,{:=}\,\deg f$.}}
\ms
This follows from Corollary~\hl{C2}{2}, and Proposition~\hl{P1.4}{1.4} below using {\it symmetries\1} of the $\mu'_k$ and $\delta''_k:=\mu''_k\mi\nu_{k+d}$ with center $3d/2$ and $d$ respectively (see \cite[(16)]{wh}) together with \cite[Thm.\,4.1]{DiPo} for (c)\,$\Leftrightarrow$\,(e), see \hl{2.5}{2.5} below.
\sk
Proposition~1 does not hold for $n\gess 4$, and there are many counterexamples (since the center of symmetry of the $\mu'_k$ is at least $2d$), see for instance \cite[Ex.\,5.6]{wh}, which shows the failure of (d)\,$\Rightarrow$\,(c) with $n\eq 4$, $d\eq 6$, $M'_6\eq M'_{18} \eq 0$, $M'_{12}\nes0$. We can also set $m\eq d{-}2\eq 3,4,5\dots$ with $n\eq 4$ in the example written just after Corollary~\hl{C3}{3} below using \cite[A.3]{wh}, where $\gamma_d\less\mu_Z$. As for the failure of (c)\,$\Rightarrow$\,(e), consider $f\eq x^d\pl y^{d-1}w\pl z^{d-1}w$ explained before Proposition~\hl{P1}{1} with $d\gess 4$, where the roots of $b_h(s)$ with $h\,{:=}\,f|_{w=1}$ (and hence those of $b_f(s)$) are {\it not\1} contained in $(-2,0)$, although $-1$ is the unique integral root of $b_f(s)$, see Remarks~\hl{R1.4b}{1.4b} and \hl{R2.5}{2.5} below. (For $n\eq 4$, there is a counterexample to \cite[Thm.\,4.1]{DiPo}, but the assumption of Corollary~\hl{C2}{2} is not satisfied.) Note also that under the first two hypotheses of Corollary~\hl{C2}{2} with $n\eq3$, we have $H^0_{\mm}\bl(\OO_{X,0}/(\dd f)\br)\eq M'\eq0$ if and only if $D\sst X$ is a free divisor, see \cite{DiSt1}. (This also fails for $n\gess 4$.)
\sk
By \cite[Thm.\,5.2]{DS2}, the pole order spectral sequence {\it never\1} degenerates at $E_2$ if the first two assumptions in Corollary~\hl{C2}{2} are satisfied, but not the last one (that is, some of the isolated singularities is not weighted homogeneous). However, this does not immediately imply the non-injectivity of (\hl{7}{7}), since we need the non-degeneration exactly at the degree $d$ part as is shown by Theorem~\hl{T3}{3}\,(i) (see also Remark~\hl{R2.4}{2.4} below).
\sk
Corollary~\hl{C2}{2} combined with Proposition~\hl{P1.4}{1.4} below, \cite[Thm.\,3]{wh}, \cite[Cor.\,1]{DS2} gives the following.
\par\htt{C3}{}\msn
{\bf Corollary~3.} {\it In the notation and assumption of Corollary~{\rm\hl{C2}{2}}, let $\mu_Z$ be the sum of Milnor $($or Tjurina$)$ numbers of isolated singularities of $Z$. Assume $\binom{d-1}{n-1}\,{>}\,\mu_Z$ with $d\,{:=}\,\deg f$, or more generally, $\gamma_{dn'}\,{>}\,\mu_Z$ with $n'\,{:=}\,[n/2]$ in the notation of {\rm(\hl{2.4.3}{2.4.3})} below. Then the morphism~{\rm(\hl{7}{7})} is injective for any $j\ins\Z$, but not surjective for $j\eq n{-}1$ or $n$.}
\ms
The hypothesis is satisfied in the case $Z$ has only one singular point which is a {\it homogeneous\1} ordinary $m$-ple point (for instance, $f\eq\msum_{i=1}^{n-1}\,x_i^m(x_i^{d-m}\pl x_n^{d-m})$) with $\tbinom{d-1}{n-1}\,{>\,}(m{-}1)^{n-1}$ and $n\gess 3$. Note that $\tbinom{d-1}{n-1}\eq\gamma_d\less\gamma_{dn'}$ in the notation of (\hl{2.4.3}{2.4.3}) below, see for instance \cite[(4.11.1)]{wh}.
\sk
Combining Theorem~\hl{T3}{3}\,(ii) with Lemma~\hl{L1.6}{1.6} and Remark~\hl{R1.6a}{1.6a}--\hl{R1.6b}{b} below, we can get the following.
\par\htt{C4}{}\msn
{\bf Corollary~4.} {\it Assume that $D$ is everywhere positively weighted homogeneous, and the BS polynomial $b_f(s)$ of a local defining function $f$ of $D$ has an integral root which is strictly smaller than $-1$. Then the morphism $(\hl{7}{7})$ is not an isomorphism for some $j\ins[1,n]$ $($which does not necessarily coincide with $n)$.}
\ms
In the isolated singularity case, this is essentially known by \cite{HM}, since the roots of $b_f(s)$ are described by the Jacobian algebra in the weighted homogeneous isolated singularity case, see Remark~\hl{R1.4b}{1.4b} below. So Corollary~\hl{C4}{4} may be viewed as a partial generalization of \cite{HM}. Note that the last hypothesis on integral roots cannot be satisfied in the case of free divisors or hyperplane arrangements by \cite{Na2}, \cite{Wa} (or \cite{bha}), and a quasi-isomorphism holds in these cases, see \cite{CNM}, \cite{Ba} (and \S\hl{S3}{3} below).
\sk
We thank the referee for useful comments to improve the paper.
\sk
In Section~\hl{S1}{1}, we review some basics of Brieskorn modules, Gauss-Manin systems, pole order filtration, and logarithmic complexes. Some consequences of the key Proposition~\hl{P1.4}{1.4} in the positively weighted homogeneous case are explained as well. In Section~\hl{S2}{2}, we prove the main theorems applying the assertions in the previous section. In Section~\hl{S3}{3}, we give a simple proof of a stronger version of the comparison theorem for hyperplane arrangements as an immediate corollary of \cite{Sc}, \cite{DeSi}. In \hl{Ap}{Appendix}, the annihilator of $f^{\al}$ is studied.
\bs\bs\htt{S1}{}
\vbox{\centerline{\bf 1. Preliminaries}
\bsn
In this section, we review some basics of Brieskorn modules, Gauss-Manin systems, pole order filtration, and logarithmic complexes. Some consequences of the key Proposition~\hl{P1.4}{1.4} in the positively weighted homogeneous case are explained as well.}
\par\htt{1.1}{}\msn
{\bf 1.1.~Gauss-Manin systems.} Let $f$ be a holomorphic function on a complex manifold $X$ of dimension $n\gess 2$. Put $D:=f^{-1}(0)\sst X$. Let $i_f:X\into X{\times}\C$ be the graph embedding by $f$ with $t$ the coordinate of $\C$. Set
$$\B_f:=(i_f)_*^{\D}\OO_X=\OO_X[\dd_t]\delta(t{-}f),$$
\par\nin where $(i_f)_*^{\D}$ is the direct image as $\D$-module (and the sheaf-theoretic direct image by $i_f$ is omitted to simplify the notation). The last term is a free module over $\OO_X[\dd_t]$ generated by $\delta(t{-}f)$, and the actions of $\dd_{x_i}$, $t$ on $\delta(t{-}f)$ are given by
\htt{1.1.1}{}
$$\dd_{x_i}\delta(t{-}f)=-(\dd_{x_i}f)\dd_t\delta(t{-}f),\q t\delta(t{-}f)=f\delta(t{-}f).
\leqno(1.1.1)$$
\par\nin \sk
The {\it Gauss-Manin systems\1} are defined by
\htt{1.1.2}{}
$$\G_{f,p}^j:=\Hc^j\K_{f,p}^{\ssb}\q\h{with}\q\K_f^{\ssb}:=\DR_{X\times\C/\C}(\B_f)[-n]\q\q(j\ins\Z,\,p\ins D),
\leqno(1.1.2)$$
\par\nin where $\DR_{X\times\C/\C}$ is the (shifted) relative de Rham complex so that $\K_f^j\eq0$ ($j\,{\notin}\,[0,n]$). By (\hl{1.1.1}{1.1.1}) the differential of $\K_f^{\ssb}$ is given by
\htt{1.1.3}{}
$$\ddd\bl(\eta\1\dd_t^k\delta(t{-}f)\br)=(\ddd\eta)\1\dd_t^k\delta(t{-}f)-(\ddd f{\wedge}\eta)\1\dd_t^{k+1}\delta(t{-}f)\q(k\ins\N),
\leqno(1.1.3)$$
\par\nin and $\K_f^{\ssb}$ is essentially the double complex associated with $\ddd$ and $\ddd f\wedge$. The $\G_{f,p}^j$ are regular holonomic $\D_{\C,0}$-modules (with $\D_{\C,0}\eq\C\{t\}\langle\dd_t\rangle$), and correspond to the Milnor cohomology groups of $f$ at $p$ using the de~Rham functor $\DR_{\C}$. They are finite free over $\C\{\!\{\dd_t^{-1}\}\!\}[\dd_t]$ for $j\ne 1$ (since the Milnor fibers are contractible). We have the isomorphisms
\htt{1.1.4}{}
$$\Gr_V^{\al}\G_{f,p}^j=H^{j-1}(F_{\!f,p},\C)_{\la}\q\q(j\ins\Z,\,\al\ins\Q,\,\la\eq e^{-2\pi i\al}).
\leqno(1.1.4)$$
\par\nin Here $F_{\!f,p}$ is the Milnor fiber of $f$ around $p$, $E_{\la}$ denotes the $\la$-eigenspace for a vector space $E$ endowed with a monodromy action, and $V$ is the filtration of Kashiwara and Malgrange, see \cite{BS}, \cite{DS2}, \cite{nwh}.
\par\htt{1.2}{}\msn
{\bf 1.2.~Brieskorn modules.} In the above notation, set
$$\A_f^j:={\rm Ker}(\ddd f\wedge:\Om_X^j\to\Om_X^{j+1}).$$
\par\nin Then $(\A_f^{\ssb},\ddd)$ is a subcomplex of $\K_f^{\ssb}$, and we have the canonical morphisms
\htt{1.2.1}{}
$$\iota_{f,p}^j:\Hc_{f,p}^j:=H^j(\A_{f,p}^{\ssb},\ddd)\to\G_{f,p}^j\q\q(j\ins\Z,\,p\ins D).
\leqno(1.2.1)$$
\par\nin The $\Hc_{f,p}^j$ are called the {\it Brieskorn modules\1} of $f$ at $p$, see \cite{Br}, \cite{BS}. They are modules over $\C\{t\}$ and also over $\C\{\!\{\dd_t^{-1}\}\!\}$. By \cite[Thm.\,1]{BS}, ${\rm Ker}\,\iota_{f,p}^j$ coincides with the $t$-torsion and also with the $\dd_t^{-1}$-torsion, and ${\rm Im}\,\iota_{f,p}^j$ is a finite free module over $\C\{t\}$ and also over $\C\{\!\{\dd_t^{-1}\}\!\}$. Its rank coincides with the dimension of the Milnor fiber cohomology, and it generates $\G_{f,p}^j$ over $\C[\dd_t]$ (more precisely, $\G_{f,p}^j\eq\mcup_k\,\dd_t^k({\rm Im}\,\iota_{f,p}^j)$), and is contained in $V^{>0}\G_{f,p}^j$ with $V$ the filtration of Kashiwara and Malgrange (which is shifted by 1 in \cite{BS}).
\sk
The action of $\dd_t^{-1}$ on $\Hc_{f,p}^j$ is given by
\htt{1.2.2}{}
$$\dd_t^{-1}[\eta]=[\ddd f{\wedge}\eta']\q\q\h{with}\q\q\ddd\eta'\eq\eta.
\leqno(1.2.2)$$
\par\nin This is compatible with (\hl{1.1.3}{1.1.3}) for $k\eq 0$.
\par\htt{R1.2}{}\msn
{\bf Remark~1.2.} Assume $f$ is positively weighted homogeneous with weights $w_1,\dots,w_n$ around a point $p\ins D$ as in the introduction. Then $\Hc_{f,p}^j$, $\G_{f,p}^j$ are completions of graded modules (with degrees in $\Q$) so that
\htt{1.2.3}{}
$$\Hc_{f,p}^j=\widehat{\mopl}_{\al\ins\Q}\,\bl(\Hc_{f,p}^j\br)_{\al},\q\G_{f,p}^j=\widehat{\mopl}_{\al\ins\Q}\,\bl(\G_{f,p}^j\br)_{\al}.
\leqno(1.2.3)$$
\par\nin Here $\bl(\Hc_{f,p}^j\br)_{\al}$, $\bl(\G_{f,p}^j\br)_{\al}$ denote the degree $\al$ part on which the Lie derivation $L_{\xi}$ is given by multiplication by $\al$ (with $\xi$ as in \hl{1.4}{1.4} below). Combining (\hl{1.2.2}{1.2.2}) with (\hl{1.4.2}{1.4.2}--\hl{1.4.3}{3}) below, we see that $\Hc_{f,p}^j$ is stable by the action of $\dd_tt$, and the latter is given by multiplication by $\al$ on the degree $\al$ part (hence this holds also for $\G_{f,p}^j$). We then get the canonical isomorphisms
\htt{1.2.4}{}
$$\Gr_V^{\al}\Hc_{f,p}^j=\bl(\Hc_{f,p}^j\br)_{\al},\q\Gr_V^{\al}\G_{f,p}^j=\bl(\G_{f,p}^j\br)_{\al}\q(\al\ins\Q),
\leqno(1.2.4)$$
\par\nin where $V$ is the filtration of Kashiwara and Malgrange such that the action of $\dd_tt\mi\al$ is nilpotent on $\Gr_V^{\al}$.
\sk
Via (\hl{1.1.4}{1.1.4}) we have the canonical isomorphisms
\htt{1.2.5}{}
$$\aligned&{\rm Im}\bl(\Gr_V^{\al}\Hc_{f,p}^j\,{\to}\,\Gr_V^{\al}\G_{f,p}^j\br)=P^kH^{j-1}(F_{\!f,p},\C)_{\la}\\ &\q\q\q(\al\ins\Q,\,[j\mi\al]\eq k,\,\la\eq e^{-2\pi i\al}),\endaligned
\leqno(1.2.5)$$
\par\nin with $P$ the pole order filtration explained in \hl{1.3}{1.3} just below. (This isomorphism can be shown using acyclicity of the complex $(\Om_X^{\ssb},\ddd)$.) In the isolated singularity case, this is quite well known (see for instance \cite{SS}), where the pole order filtration $P$ coincides with the Hodge filtration $F$.
\par\htt{1.3}{}\msn
{\bf 1.3.~Pole order spectral sequences.} We have the pole order filtration $P$ on $\K_f^{\ssb}$ defined by
$$P_k\K_f^j:=F_{k+j}\B_f{\otimes}_{\OO_X}\Om_X^j\q(k\ins\Z,\,j\ins[0,n]),$$
\par\nin where the filtration $F$ on $\B_f$ is by the order of $\dd_t$.
The $\Gr^P_k\K_f^j$ are truncated Koszul complexes for the action $\ddd f\wedge$ on $\Om_X^{\ssb}$. (In the isolated singularity case, it gives the Hodge filtration $F$, but this does not hold in the non-isolated singularity case. We have only the inclusion $F\sst P$, see for instance \cite[Prop.\,4.4]{bCM}.) We then get the pole order spectral sequence, which is essentially the spectral sequence for a double complex with differential given by $\ddd$ and $\ddd f\wedge$, see \cite{nwh}.
\sk
Assume $f$ is positively weighted homogeneous with weights $w_1,\dots,w_n$ around a point $p\ins D$ as in Remark~\hl{R1.2}{1.2} above. The spectral sequence is compatible with the weighted grading. We get the induced pole order spectral sequence on each degree $\al$ part. By (\hl{1.1.4}{1.1.4}) and (\hl{1.2.4}{1.2.4}), this spectral sequence defines the pole order filtration $P$ on $H^j(F_{\!f,p},\C)_{\la}$ via the isomorphism (\hl{1.1.4}{1.1.4}) with $\al\in(-1,0]$, where $P^k\eq P_{-k}$. We have the following.
\par\htt{T1.3}{}\msn
{\bf Theorem~1.3} (\cite[Thm.\,2]{bCM}). {\it Assume $-\al{-}k$ is not a root of $b_{f,q}(s)$ for any $k\ins\N$ and $q\ne p$ sufficiently near $q$. Then $-\al$ is a root of $b_{f,p}(s)$ if and only if}
\htt{1.3.1}{}
$$\Gr_P^jH^{n-1}(F_{\!f,p},\C)_{\la}\ne 0\q(j\eq[n\mi\al],\,\la\eq e^{-2\pi i\al}).
\leqno(1.3.1)$$
\par\nin \sk
Here $b_{f,p}(s)$ is the local BS polynomial of $f$ at $p$, and the filtration $\widetilde{P}$ coincides with $P$ in the weighted homogeneous, see a remark after \cite[(4.1.6)]{bCM}.
\par\htt{R1.3a}{}\msn
{\bf Remark~1.3a.} Assume $f$ is a homogeneous polynomial, and the projective hypersurface $Z$ defined by $f$ has only weighted homogeneous isolated singularities so that the pole order spectral sequence degenerates at $E_2$ (see \cite{wh}). Setting
$$\aligned M:=H^n_{\ddd f\wedge}(\Om^{\ssb}),&\q M^{(2)}:=H^n_{\ddd}(H^{\ssb}_{\ddd f\wedge}(\Om^{\ssb})),\\
N:=H^{n-1}_{\ddd f\wedge}(\Om^{\ssb})(-d),&\q N^{(2)}:=H^{n-1}_{\ddd}(H^{\ssb}_{\ddd f\wedge}(\Om^{\ssb}))(-d),\endaligned$$
\par\nin we have the isomorphisms for $k\in[1,d]$, $j\ins\N$, and $\la\,{:=}\, e^{-2\pi ik/d}$\,:
\htt{1.3.2}{}
$$\aligned M^{(2)}_{k+jd}&=\Gr_P^{n-1-j}H^{n-1}(F_{\!f,0},\C)_{\la},\\ N^{(2)}_{k+jd}&=\Gr_P^{n-1-j}H^{n-2}(F_{\!f,0},\C)_{\la}.\endaligned
\leqno(1.3.2)$$
\par\nin Here $\Om^{\ssb}$ denotes the complex of algebraic differential forms on $X\eq\C^n$ (which is identified with the graded quotients of the $\mm$-adic filtration on $\Om_{X,0}^{\ssb}$), and $H^{\ssb}_{\ddd f\wedge}$ means that we take the cohomology using the differential $\ddd f\wedge$ (similarly for $H^{\ssb}_{\ddd}$ with $\ddd f\wedge$ replaced by $\ddd$), see for instance \cite[(3)]{nwh}. Under the hypothesis of Theorem~\hl{T1.3}{1.3}, the first isomorphism in (\hl{1.3.2}{1.3.2}) means that $-k/d$ is a root of $b_f(s)$ if and only if $M^{(2)}_k\nes0$.
\par\htt{R1.3b}{}\msn
{\bf Remark~1.3b.} The above construction can be generalized to the case where the condition on the projective hypersurface $Z$ is not satisfied so that the pole order spectral sequence does not necessarily degenerates at $E_2$. If the projective hypersurface $Z$ has only isolated singularities, there are graded subquotients $M^{(r+1)}$, $N^{(r+1)}$ of $M^{(r)}$, $N^{(r)}$ (with $M^{(1)}\eq M$, $N^{(1)}\eq N$) obtained by taking the cohomology of the $E_r$-differential $\ddd_r$ (shifting the degree by $-rd$) inductively for $r\gess 1$, see for instance the introduction of \cite{nwh}.
\par\htt{R1.3c}{}\msn
{\bf Remark~1.3c.} Assume $f$ is a homogeneous polynomial, and the projective hypersurface $Z\sst\PP^{n-1}$ has only isolated singularities. Then we have the injectivity of the {\it cospecialization\1} morphism
\htt{1.3.3}{}
$$H^{n-2}(F_{\!f,0},\C)\into\rlap{\raise-7pt\h{$\scriptstyle p_i\to\1 0$}}\lim\,\mopl_i\,H^{n-2}(F_{\!f_i,p_i},\C)^{T_i}.
\leqno(1.3.3)$$
\par\nin Here $f_i$ is the restriction of $f$ to a transversal slice $X_i$ to each irreducible component $\Sigma_i$ of the singular locus $\Sigma\,{:=}\,{\rm Sing}\,f$ with $\{p_i\}\eq\Sigma_i\cap X_i$, and $T_i$ is the horizontal monodromy around $0\ins\Sigma_i\,(\cong\C)$, see for instance \cite[Thm.\,0.1]{DS1}. We may assume that the $X_i$ are defined by $y\eq c$ (independently of $i$) with $y$ a sufficiently general linear form on $X\eq\C^n$. Then the above limit can be obtained by putting $c\tos 0$, and we can replace the limit by the nearby cycle functor $\psi_y$. The above morphism is induced by the {\it cospecialization\1} morphism $i_0^{\prime*}\to\psi_y$ applied to the nearby cycle complex $\psi_f\C_X$, where $i'_0:\{y\eq0\}\into X$ denotes the inclusion.
\par\htt{1.4}{}\msn
{\bf 1.4.~Logarithmic complexes in the weighted homogeneous case.} Assume $X\eq\C^n$, and $D\sst X$ is a divisor defined by a reduced weighted homogeneous polynomial $f$ with strictly positive weights $w_i$ as in the introduction. Set
$$\xi:=\msum_{i=1}^n\,w_i\1x_i\dd_{x_i}.$$
\par\nin \vskip-2mm\nin
Then
\htt{1.4.1}{}
$$\iota_{\xi}(\ddd f)\eq\xi(f)\eq f,
\leqno(1.4.1)$$
\par\nin with $\iota_{\xi}$ the {\it interior product.} There are well-known relations
\htt{1.4.2}{}
$$\iota_{\xi}\ssc\ddd+\ddd\ssc \iota_{\xi}=L_{\xi},
\leqno(1.4.2)$$
\par\nin \vskip-7mm
\htt{1.4.3}{}
$$\iota_{\xi}\ssc(\ddd f\!/\!f\wedge)+(\ddd f\!/\!f\wedge)\ssc \iota_{\xi}={\rm id},
\leqno(1.4.3)$$
\par\nin where $L_{\xi}$ denotes the {\it Lie derivation.} (The last equality follows from the Leibniz rule.)
\sk
The {\it logarithmic forms\1} $\Om_X^j(\log D)$ are defined by the conditions: $f\eta\ins\Om_X^j$ and $f\ddd\eta\ins\Om_X^{j+1}$ for $\eta\ins\Om_X^j(*D)$, see \cite{SaK}. The last condition is equivalent to that $\ddd f{\wedge}\eta\ins\Om_X^{j+1}$ assuming the first. Using (\hl{1.4.3}{1.4.3}), it is easy to see the following (see also \cite{HM}).
\par\htt{L1.4}{}\msn
{\bf Lemma~1.4.} {\it Let $\A_f^j$ be as in {\rm\hl{1.2}{1.2}}. There are decompositions
\htt{1.4.4}{}
$$\Om_X^j(\log D)=\A_f^jf^{-1}\oplus \iota_{\xi}\A_f^{j+1}\!f^{-1}\q\q(j\ins\Z),
\leqno(1.4.4)$$
\par\nin together with the isomorphisms}
\htt{1.4.5}{}
$$\iota_{\xi}:\A_f^{j+1}\!f^{-1}\,\simto\,\iota_{\xi}\A_f^{j+1}\!f^{-1}\q\q(j\ins\Z).
\leqno(1.4.5)$$
\par\nin \msn
{\it Proof.} By (\hl{1.4.3}{1.4.3}) the identity on the complex $(\Om_X^{\ssb}(\log D),\ddd f\!/\!f\wedge)$ is homotopic to 0. Hence the complex is acyclic, and a splitting of complex is given by $\iota_{\xi}$ using (\hl{1.4.3}{1.4.3}). Lemma~\hl{L1.4}{1.4} thus follows.
\par\htt{P1.4}{}\msn\vbox{\nin
{\bf Proposition~1.4.} {\it Let $\bl(\Om_X^{\ssb}(\log D)f^{-r},\ddd\br)$ be the logarithmic complex multiplied by $f^{-r}$ for $r\ins\N$. This complex is isomorphic to the mapping cone
$$C\bl(L_{\xi}:(\A_f^{\ssb}f^{-r-1},\ddd)\to(\A_f^{\ssb}f^{-r-1},\ddd)\br),$$
\par\nin and in the notation of Remark~{\rm\hl{R1.2}{1.2}} we have the isomorphisms for $j\ins\Z\,{:}$}
\htt{1.4.6}{}
$$H^j\bl(\Om_{X,0}^{\ssb}(\log D)f^{-r},\ddd\br)=\bl(\Hc_{f,0}^j\br)_{r+1}\oplus\bl(\Hc_{f,0}^{j+1}\br)_{r+1}.
\leqno(1.4.6)$$}
\skn
{\it Proof.} It is enough to show the first assertion (using Remark~\hl{R1.2}{1.2}). Indeed, the action of $L_{\xi}$ on the degree $\al$ part is given by multiplication by $\al$. So we may restrict to the degree 0 part (since the other part is acyclic). Here the mapping cone is associated with the zero map, hence it is a direct sum of two complexes. Note also that the differential $\ddd$ on $\A_{f,p}^{\ssb}$ commutes with multiplication by $f\1^r$ ($r\ins\N$).
\sk
To show the first assertion, we first see that $\bl(\Om_X^{\ssb}(\log D)f^{-r},\ddd\br)$ contains $(\A_f^{\ssb}f^{-r-1},\ddd)$ as a subcomplex by Lemma~\hl{L1.4}{1.4}, since $\ddd$ and $\ddd f\wedge$ anti-commute.
\sk
On the other hand, the restriction of $\ddd$ to $\iota_{\xi}\A_f^{j+1}\!f^{-r-1}$ is the sum of
$$\ddd':\iota_{\xi}\A_f^{j+1}\!f^{-r-1}\to\A_f^{j+1}\!f^{-r-1}\q\h{and}\q\ddd'':\iota_{\xi}\A_f^{j+1}\!f^{-r-1}\to \iota_{\xi}\A_f^{j+2}\!f^{-r-1}.$$
\par\nin Using (\hl{1.4.2}{1.4.2}) and the diagram below, we see that $\ddd'$, $\ddd''$ are identified respectively with $L_{\xi}$ and the restriction of $\ddd$ to $\A_f^j\!f^{-r-1}$ up to sign via the isomorphism (\hl{1.4.5}{1.4.5}).
\htt{1.4.7}{}
$$\begin{array}{ccccc}\buildrel{\!\ddd}\over\to&\A_f^{j+1}\!f^{-r-1}&\buildrel{\!\ddd}\over\to&\A_f^{j+2}\!f^{-r-1}&\buildrel{\!\ddd}\over\to\\&\raise2pt\h{$\scriptstyle \iota_{\xi}$}\!\downarrow\!\uparrow\!\raise1pt\h{$\scriptstyle\ddd'$}&\raise5mm\h{}\raise-2mm\h{}&\raise2pt\h{$\scriptstyle \iota_{\xi}$}\!\downarrow\!\uparrow\!\raise1pt\h{$\scriptstyle\ddd'$}\\ \buildrel{\ddd''}\over\to&\iota_{\xi}\A_f^{j+1}\!f^{-r-1}&\buildrel{\ddd''}\over\to&\iota_{\xi}\A_f^{j+2}\!f^{-r-1}&\buildrel{\ddd''}\over\to\end{array}
\leqno(1.4.7)$$
\par\nin So Proposition~\hl{P1.4}{1.4} follows.
\par\htt{R1.4a}{}\msn
{\bf Remark~1.4a.} If $f$ is a weighted homogeneous polynomial with an {\it isolated\1} singularity at 0 and $n\gess 3$, then Proposition~\hl{P1.4}{1.4} together with (\hl{1.2.4}{1.2.4}--\hl{1.2.5}{5}) and Remark~\hl{R1.4b}{1.4b} just below implies that the morphism (\hl{2}{2}) is a quasi-isomorphism if and only if the {\it unipotent\1} monodromy part of the vanishing cohomology vanishes, or equivalently, there is no {\it integral\1} spectral number (or no {\it integral\1} root of the {\it reduced\1} BS polynomial $\bt_f(s)\,{:=}\,b_f(s)/(s{+}1)$). In the weighted homogeneous isolated singularity case with $n\eq2$, the morphism (\hl{2}{2}) is always a quasi-isomorphism using Proposition~\hl{P1.4}{1.4} (where the link may be disconnected). These imply another proof of Theorem in \cite{HM}.
\par\htt{R1.4b}{}\msn
{\bf Remark~1.4b.} Assume $f$ is a weighted homogeneous polynomial with an isolated singularity at 0 and with weights $w_i$. The spectrum ${\rm Sp}_f(t)\eq\msum_k\,t\1^{\al_{f,k}}$ has a {\it symmetry}
\htt{1.4.8}{}
$${\rm Sp}_f(t)={\rm Sp}_f(t^{-1})\1t^n,
\leqno(1.4.8)$$
\par\nin which holds in the general hypersurface isolated singularity case, and it coincides with the Poincar\'e series of the Jacobian ring $\C\{x\}/(\dd f)$ shifted by $\msum_i\,w_i$. More precisely, we have
\htt{1.4.9}{}
$${\rm Sp}_f(t)=\prod_{i=1}^n\1\frac{t-t^{w_i}}{t^{w_i}\mi 1},
\leqno(1.4.9)$$
\par\nin see \cite{St1}, \cite{St2} (and also \cite{JKSY}).
\sk
The spectral numbers $\al_{f,k}$ coincide with the roots of the reduced BS polynomial $\bt_f(s)$ up to sign (forgetting multiplicities). This follows by combining \cite{Va} (or \cite{SS}) with \cite{Ma}. Note that there is no {\it integral\1} spectral number (that is, the {\it unipotent\1} monodromy part of the vanishing cohomology vanishes) if and only if $D$ is a $\Q$-{\it homology manifold.} This follows from the Wang sequence, see \cite{Mi}.
\par\htt{R1.4c}{}\msn
{\bf Remark~1.4c.} The inductive limit of the isomorphism (\hl{1.4.6}{1.4.6}) over $r\ins\N$ is closely related to the assertion that $X\stm D$ is a $\C^*$-bundle over $\PP^{n-1}\stm Z$ in the $f$ homogeneous polynomial case (using for instance \cite[\S1.3]{BuSa}).
\ms
Lemma~\hl{L1.4}{1.4} has the following.
\par\htt{C1.4}{}\msn
{\bf Corollary~1.4.} {\it Assume $X\eq\C^n$ with $D\sst X$ defined by a positively weighted homogeneous reduced polynomial $f$ as in {\rm\hl{1.4}{1.4}}, and $\dim{\rm Sing}\,D\less 1$ $($for instance, $n\eq 3)$. Then $D$ is tame, that is, the logarithmic differential forms $\Om_X^j(\log D)$ have at most projective dimension $j$ for any $j\ins\Z$.}
\msn
{\it Proof.} By Lemma~\hl{L1.4}{1.4}, it is enough to show that
\htt{1.4.10}{}
$${\rm pd}\1_{\OO_X}\A_f^j<j\q\h{if}\q j\less n{-}1,
\leqno(1.4.10)$$
\par\nin since this is clear for $j\eq n$ by definition (that is, $\A_f^n\eq\Om_X^n(D)$). Here we use algebraic coherent sheaves, and ${\rm pd}\1_{\OO_X}$ denotes the projective dimension over $\OO_X$. It is well known that we have {\it acyclicity\1} of the Koszul complex
\htt{1.4.11}{}
$$\Hc^j(\Om_X^{\ssb},\ddd f\wedge)\eq0\q\h{if}\q j\,{<}\,n{-}1,
\leqno(1.4.11)$$
\par\nin since $\dim{\rm Sing}\,D\less 1$, see for instance \cite[Prop.\,2]{JKSY2}. This implies the assertion (\hl{1.4.10}{1.4.10}) for $j\,{<}\,n{-}1$, since $\A_f^j\eq{\rm Ker}\,\ddd f{\wedge}\eq{\rm Im}\,\ddd f{\wedge}\,(\subset\Om_X^j)$ for such $j$.
\sk
For $j\eq n{-}1$, we have the exact sequence
\htt{1.4.12}{}
$$0\to\A_f^{n-1}\to\Om_X^{n-1}\buildrel{\!\!\!\ddd f\wedge}\over\longrightarrow\Om_X^n\tos M^{\sim}\to 0,
\leqno(1.4.12)$$
\par\nin where $M^{\sim}$ is the coherent sheaf associated to the graded $\C[x_1,\dots,x_n]$-module $M$ studied in \cite{DS2}, \cite{wh}, \cite{nwh} (and is defined by the exact sequence). We then get that
\htt{1.4.13}{}
$$\aligned{\rm pd}\1_{\OO_X}\A_f^{n-1}&={\rm pd}\1_{\OO_X}\I^{(n)}\mi 1\\ &={\rm pd}\1_{\OO_X}M^{\sim}\mi 2\,\les\,n{-}2,\endaligned
\leqno(1.4.13)$$since $X$ is smooth, where we denote by $\I^{(n)}\sst\Om_X^n$ the image of $\ddd f\wedge$. This finishes the proof of Corollary~\hl{C1.4}{1.4}.
\par\htt{1.5}{}\msn\vbox{\nin
{\bf 1.5.~Twisted logarithmic complexes.} Assume $D\sst X$ is everywhere positively weighted homogenous. In the notation of the introduction, set
$$\om_p:=\msum_k\,\al_{p,k}\1\om_{p,k}\q\q\h{with}\q\q\om_{p,k}:=\ddd f_{p,k}/\!f_{p,k}.$$
\par\nin We have the following.
\par\htt{P1.5}{}\msn
{\bf Proposition~1.5.} {\it The complex $\bl(\Om_{X,p}^{\ssb}(\log D)f^{-r},\ddd\pl\om_p\wedge\br)$ for $p\ins D$ and $r\ins\N$ is isomorphic to the mapping cone
\htt{1.5.1}{}
$$C\bl(L_{\xi_p}\pl\alt_p:(\A_{f,p}^{\ssb}f^{-r-1},\ddd\pl\om_p\wedge)\to(\A_{f,p}^{\ssb}f^{-r-1},\ddd\pl\om_p\wedge)\br).
\leqno(1.5.1)$$
\par\nin Here $\xi_p$ is the vector field associated with a positively weighted homogeneous polynomial $f_p$ as in {\rm\hl{1.4}{1.4}}, and $\alt_p$ is defined in $(\hl{4}{4})$.}}
\msn
{\it Proof.} The argument is essentially the same as in the proof of Proposition~\hl{P1.4}{1.4}. We can calculate the action of $\om_{p,k}\wedge$ in a similar way to the case of the differential $\ddd$ using (\hl{1.4.3}{1.4.3}) instead of (\hl{1.4.2}{1.4.2}). Here $f$ is replaced by $f_{p,k}$, and the right-hand side of (\hl{1.4.3}{1.4.3}) becomes $d_{p,k}$. Proposition~\hl{P1.5}{1.5} then follows.
\par\htt{C1.5}{}\msn
{\bf Corollary~1.5.} {\it Assume $D$ is everywhere positively weighted homogenous. Let $\Cc^{\prime\ssb}$ be the mapping cone of the comparison morphism {\rm(\hl{2}{2})\rm}. For $p\ins D$, $j\ins\Z$, we have}
\htt{1.5.2}{}
$$\dim_{\C}\Hc^{j-1}\Cc^{\prime\ssb}_p\ges\dim_{\C}\Hc^j\Cc^{\prime\ssb}_p\q\h{if}\q\Hc^{j+1}\Cc^{\prime\ssb}_p\eq0.
\leqno(1.5.2)$$
\par\nin \msn
{\it Proof.} Proposition~\hl{P1.5}{1.5} implies that $\Cc^{\prime\ssb}_p$ is isomorphic to the mapping cone of the action of $L_{\xi_p}{+}\1\alt_p$ on the complex
$$\K_p^{\ssb}:=\bl(\rlap{\raise-10pt\h{$\,\,\,\scriptstyle r$}}\indlim\,\A^{\ssb}_{f,p}f^{-r-1}/\A^{\ssb}_{f,p}f^{-1},\ddd\pl\om_p{\wedge}\br).$$
\par\nin Here the action of $L_{\xi_p}$ on the modified degree $\beta$ part of $\A^{\ssb}_{f,p}f^{-r}$ is given by multiplication by $\beta\ins\C$ by the definition of modified degree, see Remark~\hl{R1.2}{1.2} for the untwisted case. Taking the kernel or cokernel of the action of $L_{\xi_p}{+}\1\alt_p$ is then the same as the restriction to the modified degree $-\alt_p$ part. (Note that the differential $\ddd\pl\om_p{\wedge}$ preserves the modified degree.) Using the long exact sequence associated with the mapping cone of the action of $L_{\xi_p}{+}\1\alt_p$ on $\K_p^{\ssb}$, we get the inequality (\hl{1.5.2}{1.5.2}). This finishes the proof of Corollary~\hl{C1.5}{1.5}.
\par\htt{1.6}{}\msn
{\bf 1.6.~Logarithmic stratification.} For a reduced divisor $D$ on a complex manifold $X$, we denote by $\Theta_X(-\log D)$ the sheaf of {\it logarithmic vector fields.} By definition a holomorphic vector field $\xi\ins\Theta_X$ belongs to $\Theta_X(-\log D)$ if and only if $\xi(f)\sst(f)$ with $f$ a holomorphic function defining locally $D\sst X$, see \cite{SaK}.
\sk
For $p\in D$, let $E_p\sst T_pX$ be the image of $\Theta_{X,p}(\log D)$ in the tangent space $T_pX$ which is identified with $\Theta_{X,p}{\otimes}_{\OO_{X,p}}\OO_{X,p}/\mm_{X,p}$, where $\mm_{X,p}\sst\OO_{X,p}$ is the maximal ideal. Set
$$D^{(k)}:=\{p\in D\mid\dim E_p\eq k\}\q(k\gess 0).$$
\par\nin Taking local generators of $\Theta_X(-\log D)$, we get a matrix with coefficients in $\OO_X$, and the $D^{(k)}$ are determined by the rank of this matrix at each $p$. This implies that the union $D^{(\les j)}:=\mcup_{k\les j}\,D^{(k)}$ is a closed analytic subset for any $j\ins\N$.
\par\htt{D1.3}{}\msn
{\bf Definition~1.6.} We say that a divisor $D\sst X$ has the {\it locally finite logarithmic stratification,} if $\dim D^{(k)}\less k$ for any $k\gess0$. (This is equivalent to that the ``logarithmic stratification" in the sense of \cite{SaK} is everywhere locally finite.)
\par\htt{R1.6a}{}\msn
{\bf Remark~1.6a.} Assume $D$ has the locally finite logarithmic stratification. Then, choosing a basis $(v_1,\dots,v_k)$ of $E_p$ with $k:=\dim E_p$, there are inductively local analytic isomorphisms for $j\in[1,k]$\,:
\htt{1.6.1}{}
$$\aligned&(X,p)\cong(S^{(j)}_p,p){\times}(\Delta^j,0)\q\h{inducing}\\&(D,p)\cong(S^{(j)}_p\caps D,p){\times}(\Delta^j,0),\endaligned
\leqno(1.6.1)$$
\par\nin such that the subspace of $T_pX$ spanned by $v_1,\dots,v_j$ is identified with the tangent space of $\{p\}{\times}\Delta^j$. Here $(S^{(j)}_p,p)\sst(X,p)$ is a submanifold of codimension $j$ such that $v_1,\dots,v_j$ form a basis of $T_pX/T_pS_p$ (and $\Delta$ is an open disk). Indeed, we can show the following isomorphisms by {\it decreasing\1} induction on $j$ using integral curves of vector fields $\xi_j$ whose images in $T_pX$ are $v_j$ ($j\in[1,k]$):
\htt{1.6.2}{}
$$\aligned&(S^{(j-1)}_p,p)\cong(S^{(j)}_p,p){\times}(\Delta,0)\q\h{inducing}\\&(S^{(j-1)}\caps D,p)\cong(S^{(j)}_p\caps D,p){\times}(\Delta,0).\endaligned
\leqno(1.6.2)$$
\par\nin Here $S^{(j-1)}_p$ is the union of integral curves of $\xi_j$ passing through a point of $S^{(j)}_p$ (with $S^{(0)}_p\eq X$). This construction is compatible with $D$ by the uniqueness of integral curves (at the smooth points of $D$). Note that $S^{(j)}_p$ ($j\ins[1,k{-}1]$) is determined by the $\xi_i$ ($i\ins[j{+}1,k]$) and $S^{(k)}_p$.
\sk
In the case $D$ has the locally finite logarithmic stratification, (\hl{1.6.1}{1.6.1}) imply that the $D^{(k)}$ are smooth and purely $k$-dimensional (unless $D^{(k)}\eq\emptyset$) for any $k\gess 0$, see also \cite{SaK}.
\par\htt{R1.6b}{}\msn
{\bf Remark~1.6b.} If $D$ is positively weighted homogeneous at $p$, we may assume that so is $S^{(k)}_p\caps D\sst S^{(k)}_p$ in Remark~\hl{R1.6a}{1.6a}. Indeed, we can define $S^{(k)}_p$ by $\mcap_{i\in I}\{x_i\eq 0\}$ for some subset $I\sst\{1,\dots,n\}$ with $x_1,\dots,x_n$ {\it weighted coordinates\1} for $f_p$. (It does not seem clear whether this can be used for another proof of Lemma~\hl{L1.6}{1.6} below, since the situation at $q\ins S^{(k)}_p\caps D^{(k)}\stm\{p\}$ seems rather unclear when $k<\dim D^{(k)}$.)
\ms
The following seems be known to specialists (see for instance \cite{CNM}).
\par\htt{L1.6}{}\msn
{\bf Lemma~1.6.} {\it Assume the divisor $D$ is everywhere positively weighted homogeneous as in the introduction. Then $D$ has the locally finite logarithmic stratification.}
\msn
{\it Proof.} We have to show that $\dim D^{(k)}\less k$ for any $k\gess0$. Take a smooth point $p\ins D^{(k)}$. We have the vector field $\xi=\msum_i\,w_i\1x_i\1\dd_{x_i}$ associated to $f_p$ so that $\xi(f_p)\eq f_p$. Using the trivialization of the tangent bundle by the coordinates $x_1,\dots,x_n$, this vector field is identified with the map
\htt{1.6.3}{}
$$(x_1,\dots,x_n)\mapsto(w_1x_1,\dots,w_nx_n).
\leqno(1.6.3)$$
\par\nin If $\dim D^{(k)}\,{>}\,k\,(=\dim E_p)$, we then see that its image cannot be contained in $E_q$ for some $q\ins D^{(k)}$ sufficiently near $p$, since the $E_q$ depend on $q$ continuously in the Grassmannian. This is, however, a contradiction, since $\xi\ins\Theta_X(-\log D)$. So Lemma~\hl{L1.6}{1.6} follows.
\ms
We have the following.
\par\htt{P1.6}{}\msn
{\bf Proposition~1.6.} {\it Assume the divisor $D$ has the locally finite logarithmic stratification $($for instance, $D$ has only isolated singularities$)$. Then the morphism~{\rm(\hl{2}{2})} is a quasi-isomorphism if and only if it is a $D$-quasi-isomorphism.}
\msn
{\it Proof.} Apply the functor $\DRl_X^{-1}$ to the morphism (\hl{2}{2}) (see a remark after (\hl{3}{3}) and \cite{ind}), and consider the mapping cone
\htt{1.6.4}{}
$$\Cc^{\ssb}:=C\bl(\DRl_X^{-1}(\M_{X,\log}^{\ssb}(L),\naal)\to\DRl_X^{-1}(\M_X^{\ssb}(L),\naal)\br).
\leqno(1.6.4)$$
\par\nin This is a bounded complex of {\it coherent\1} $\D_X$-modules. We have to show that this is acyclic if $\DR_X(\Cc^{\ssb})$ is acyclic. Set $Z:=\mcup_j\,{\rm Supp}\,\Hc^j\Cc^{\ssb}$, and assume $Z\ne\emptyset$. By Remark~\hl{R1.6a}{1.6a}, $Z$ is a union of strata of $\mathcal S$. Let $V$ be a maximal-dimensional stratum contained in $Z$. By the local analytic triviality along $V$, we may restrict to an appropriate transversal slice to $V$ as in Remark~\hl{R1.6a}{1.6a} so that the assertion is reduced to the case $Z$ is a point, denoted by $0$. The cohomology sheaves $\Hc^j\Cc^{\ssb}$ are then finite direct sums of $\B_0:=\C[\dd_{x_1},\dots,\dd_{x_n}]$, see for instance \cite{Ma}. (This is a special case of Kashiwara's equivalence, see for instance \cite{reg}.) We now get a contradiction, since $\DR_X\B_0\eq\C_{\{0\}}$. This finishes the proof of Proposition~\hl{P1.6}{1.6}.
\par\htt{1.7}{}\msn
{\bf 1.7.~Holonomicity and regular holonomicity.} Using the assertions in \hl{1.6}{1.6}, we can prove the following.
\par\htt{P1.7a}{}\msn\vbox{\nin
{\bf Proposition~1.7a.} {\it Assume the divisor $D$ has the locally finite logarithmic stratification. Then the complex $\DRl_X^{-1}(\M^{\ssb}_{X,\log}(L),\nabla^{(\alpha)})$ is holonomic, that is, it belongs to $D^b_{\rm hol}(\D_X)\1;$ more precisely, the characteristic varieties of its cohomology $\D$-modules are contained in the union of the conormal bundles of $D^{(j)}$ $(j\ins[0,n])$.}}
\msn
{\it Proof.} We argue by induction on strata. Assuming the assertion is proved on the complement of $D^{\les j}$, we have to show the assertion on the complement of $D^{\les j-1}$, that is, on a neighborhood of any point $p\ins D^{(j)}$. By the analytic triviality along $D^{(j)}$, we may assume that $j\eq 0$ restricting to a transversal slice as in Remark~\hl{R1.6a}{1.6a}. Then the assertion is clear since we have
$$T^*X\stm T^*(X\stm\{p\})=T^*_pX,$$
\par\nin and $\DRl_X^{-1}(\M^{\ssb}_{X,\log}(L),\nabla^{(\alpha)})$ is a complex of coherent $\D_X$-modules by definition. Note that the locally finite logarithmic stratification satisfies the Whitney condition (b) by the local analytic triviality along strata (see Remark~\hl{R1.6a}{1.6a}), and the union of the conormal bundles of $D^{(j)}$ $(j\ins[0,n])$ is a closed subset. This finishes the proof of Proposition~\hl{P1.7a}{1.7a}
\par\htt{C1.7a}{}\msn
{\bf Corollary~1.7a.} {\it Assume $D$ has the locally finite logarithmic stratification, and is tame, see Corollary~{\rm\hl{C1.4}{1.4}}. Then the complex $\DRl_X^{-1}(\M^{\ssb}_{X,\log}(L),\nabla^{(\alpha)})[n]$ is quasi-isomorphic to a holonomic $\D_X$-module.}
\msn
{\it Proof.} Set $\K^{\ssb}:=\DRl_X^{-1}(\M^{\ssb}_{X,\log}(L),\nabla^{(\alpha)})[n]$. This is locally quasi-isomorphic to a complex of finite free $\D_X$-modules $\F^{\ssb}$ with $\F^j\eq0$ for $j\,{\notin}\,[-n,0]$ by the tameness assumption (using a double complex consisting of free resolutions locally). Applying the functor $\DD$ which assigns the dual of $\D$-modules in the derived category and preserves holonomic $\D$-modules (see for instance \cite{ind}, \cite{reg}), this property is also satisfied for $\DD\K^{\ssb}$. (Indeed, if $\F$ is a finite free $\D_X$-module, then $\Hc^j\DD\F\eq0$ for $j\nes{-}n$ by definition.) We thus get that $\DD\Hc^j\K^{\ssb}\eq\Hc^{-j}\DD\K^{\ssb}\eq0$ for $j\,{<}\,0$ (since $\Hc^j\DD\M\eq0$ for $j\nes0$ if $\M$ is holonomic). Hence $\Hc^j\K^{\ssb}\eq0$ for $j\nes0$. This finishes the proof of Corollary~\hl{C1.7a}{1.7a}.
\par\htt{P1.7b}{}\msn
{\bf Proposition~1.7b.} {\it Assume $D$ is everywhere defined by a homogeneous polynomial locally on $X$ $($as in the hyperplane arrangement case$)$. Then the complex $\DRl_X^{-1}(\M^{\ssb}_{X,\log}(L),\nabla^{(\alpha)})$ is regular holonomic, that is, it belongs to $D^b_{\rm rh}(\D_X)$.}
\msn
{\it Proof.} We argue by induction on strata of the locally finite logarithmic stratification, see Lemma~\hl{L1.6}{1.6}. Assume the assertion is proved on the complement of $D^{\les j}$. By the same argument as in the proof of Proposition~\hl{P1.7a}{1.7a}, we may assume $j\eq0$ cutting by a transversal slice to a stratum. (If there is a homogeneous polynomial $f$ together with a vector field $\xi$ not vanishing at 0 and such that $\xi(f)\in(f)$, then we can show that $f$ is a polynomial of fewer variables changing linearly the variables if necessary.)
\sk
Let $\rho\,{:}\,(\Xt,\Dt)\tos(X,D)$ be the blow-up at the origin with $\Dt\eq\rho^{-1}(D)$. Set $E\,{:=}\,\rho^{-1}(0)$. Note that $\Xt$ is a line bundle over $E$, where $E$ is identified with the zero-section, and $\Dt$ is the union of the pull-back of $D$ with the zero-section. We can then apply the K\"unneth formula for logarithmic forms locally by trivializing the line bundle.
\sk
For a sufficiently large integer $m$, there are canonical morphisms
\htt{1.7.1}{}
$$\aligned&\bl(\M_{X,\log}^{\ssb},\naal\br)\to\rho_*\bl(\M_{\Xt,\log}^{\ssb}(mE),\naal\br)\\&\q\q\q\to\R\rho_*\bl(\M_{\Xt,\log}^{\ssb}(mE),\naal\br).\endaligned
\leqno(1.7.1)$$
\par\nin Here $\R\rho_*$ is defined by taking the canonical flasque resolution (using discontinuous sections) by Godement, and $(mE)$ means ${\otimes}_{\OO_{\Xt}}\OO_{\Xt}(mE)$. The functor $\DRl^{-1}$ (see (\hl{3}{3}) for definition) commutes with the direct image by $\rho$ (see \cite[3.3--5]{ind}), and we have
$$\DRl^{-1}\bl(\M_{\Xt,\log}^{\ssb}(mE),\naal\br)\in D^b_{\rm rh}(\D_{\Xt}),$$
\par\nin using the K\"unneth formula for logarithmic forms together with the inductive hypothesis.
\sk
Let $\Cc'{}^{\ssb}$ be the mapping cone of the composition (\hl{1.7.1}{1.7.1}). This is cohomologically supported at the origin. By an argument similar to the proof of Proposition~\hl{P1.6}{1.6} (using a special case of Kashiwara's equivalence), we can then verify that
$$\DRl^{-1}\Cc'{}^{\ssb}\in D^b_{\rm rh}(\D_X).$$
\par\nin We thus get that
$$\DRl^{-1}\bl(\M_{X,\log}^{\ssb},\naal\br)\in D^b_{\rm rh}(\D_X),$$
\par\nin since regular holonomic $\D$-modules are stable by subquotients and extensions. This finishes the proof of Proposition~\hl{P1.7b}{1.7b}.
\ms
From Proposition~\hl{P1.7b}{1.7b}, Corollary~\hl{C1.7a}{1.7a}, and Lemma~\hl{L1.6}{1.6}, we get the following.
\par\htt{C1.7b}{}\msn
{\bf Corollary~1.7b.} {\it Assume $D$ is tame $($see Corollary~{\rm\hl{C1.4}{1.4})}, and is everywhere defined by a homogeneous polynomial locally on $X$ $($as in the hyperplane arrangement case$)$. Then the complex $\DRl_X^{-1}(\M^{\ssb}_{X,\log}(L),\nabla^{(\alpha)})[n]$ is quasi-isomorphic to a regular holonomic $\D_X$-module.}
\par\htt{R1.7a}{}\msn
{\bf Remark~1.7a.} The highest cohomology of $\DRl_X^{-1}(\M^{\ssb}_{X,\log}(L),\nabla^{(\alpha)})$ is the quotient of $\D_X$ divided by an ideal generated by {\it logarithmic vector fields with certain $\OO$-linear terms.} In the case $\al_k\eq\al$ ($\forall\,k$), this ideal is generated by $\Tht_{f,\al-1}$ with
$$\Tht_{f,\al}:=\bl\{\1\xi\mi\al\1\xi(f)/f\mid\xi\ins\Theta_X(-\log D)\br\}\q(\al\ins\C),$$
\par\nin using (\hl{3.1.8}{3.1.8}) below. The shift of $\al$ by $-1$ comes from the isomorphism $\Om_X^n(\log D)\eq\Om_X^n(D)$. We use the anti-involution $^*$ of $\D_X$ such that $x_i^*\eq x_i$, $\dd_{x_i}^*\eq{-}\dd_{x_i}$, and $(PQ)^*\eq Q^*P^*$ as is well known. It is interesting whether this quotient coincides with the $\D_X$-module generated by $f^{\al}$, see Theorem~\hl{T2}{2} and also Corollary~\hl{C3.2d}{3.2d}, Theorem~\hl{TA}{A} below.
\par\htt{R1.7b}{}\msn
{\bf Remark~1.7b.} It does not seem trivial to show the regularity of the highest cohomology of $\DRl_X^{-1}(\M^{\ssb}_{X,\log}(L),\nabla^{(\alpha)})$ even under the assumption of Corollary~\hl{C1.7a}{1.7a} (for instance, if it has a $\D_X$-submodule whose support has dimension 1, where we cannot apply Hartogs theorem). Note that the {\it principal symbols\1} of logarithmic vector fields do not necessarily generate the {\it reduced ideal\1} of the characteristic variety; for instance, in the case $f\eq x^d\pl y^d$ ($d\gess 3$), the logarithmic vector fields are generated by the Euler field $x\dd_x\pl y\dd_y$ and $y^{d-1}\dd_x\mi x^{d-1}\dd_y$ (since $f_x,f_y$ form a regular sequence), but $\al x\pl\beta y$ and $\al y^{d-1}\mi\beta x^{d-1}$ (with $\al,\beta\ins\C$ fixed) do not generate the maximal ideal of $\C\{x,y\}$.
\par\htt{R1.7c}{}\msn
{\bf Remark~1.7c.} We can extend Corollary~\hl{C1.7b}{1.7b} to the case where the transversal slice is defined by a positively weighted homogeneous polynomial, and has an isolated singularity. (In this case, the blowing-up at the origin can be replaced by an embedded resolution.) This includes the case where the hypothesis of Corollary~\hl{C2}{2} is satisfied (that is, $f$ is a homogeneous polynomial and $Z:=\{f\eq0\}\sst\PP^{n-1}$ has only {\it weighted homogenous\1} isolated singularities) using Corollary~\hl{C1.4}{1.4}. In the case the local system $L$ is constant, it does not seem easy to prove the assertion under the assumption of Corollary~\hl{C2}{2} using the theory of $t$-structure \cite{BBD}. Indeed, the relation between the isomorphism in Lemma~\hl{L1.4}{1.4} and the one applied to a transversal slice does not seem very clear (via the cospecialization morphism as in (\hl{1.3.3}{1.3.3})), hence it is not easy to use \cite[Prop.\,2.2]{DiSt2}.
\par\htt{R1.7d}{}\msn
{\bf Remark~1.7d.} By Corollaries~\hl{C1.4}{1.4} and \hl{C1.7a}{1.7a}, we see that under the assumption of the former, the complex $\DRl_X^{-1}(\M^{\ssb}_{X,\log}(L),\nabla^{(\alpha)})[n]$ is quasi-isomorphic to a $\D_X$-module, that is,
\htt{1.7.2}{}
$$\Hc^j\bl(\DRl_X^{-1}(\M^{\ssb}_{X,\log}(L),\nabla^{(\alpha)})\br)\eq0\q\h{if}\q j<n.
\leqno(1.7.2)$$
\par\nin Note that the {\it local weighted homogeneity\1} of the associated projective hypersurface assumed in Corollary~\hl{C2}{2} is not needed to show this quasi-isomorphism in the case $D$ is defined by a positively weighted homogeneous polynomial and $\dim{\rm Sing}\,D\less 1$.
\sk
Consider for instance the case $f\eq x^5\pl y^4z\pl x^3y^2$. Here the weighted degree 1 part of the highest Brieskorn module $(\Hc^3_{f,0})_1$ has {\it non-zero torsion,} see \cite{hil}, \cite[Ex.\,5.3]{wh}, \cite[Ex.\,5.1]{nwh}. The above assertion then implies the {\it non-injectivity\1} of the morphism (\hl{7}{7}) for $j\eq 2$ at general points of ${\rm Sing}\,D$, so the logarithmic comparison theorem fails for a non-weighted-homogeneous curve defined by $f|_{z=1}\eq0$ as is shown in \cite{CMNC} although plane curves are free divisors. (Indeed, we get a contradiction by the theory of $t$-{\it structures\1} \cite{BBD} together with Proposition~\hl{P1.4}{1.4} otherwise.) Note that $x^3y^2$ in $f$ can be replaced by $x^4y$.
\sk
It seems very interesting whether the vanishing (\hl{1.7.2}{1.7.2}) holds in the case $D$ is everywhere positively weighted homogeneous, although it seems quite difficult to prove it or to find a counterexample. Here we have also a problem of torsion of Brieskorn modules except for the ``lowest degree" part, see \cite[(4.3.6)]{DS2}.
\par\htt{R1.7e}{}\msn
{\bf Remark~1.7e.} It is proved in \cite[Ex.\,5.1]{CN0} that an example with (\hl{1.7.2}{1.7.2}) unsatisfied is given by $(xz{+}\1y)(x^4{+}\1y^5{+}\1xy^4)\eq0$. Here one may consider the homogeneous polynomial $f\eq(xz{+}\1yw)(x^4w{+}\1y^5{+}\1xy^4)$ associated with it, and calculate the {\it pole order spectral sequence\1} applying a small computer program made by using C and explained in \cite{nwh}. One then gets the non-degeneration $\nu_{2d}^{(2)}\nes\nu_{2d}^{(3)}$ in the notation of \cite{nwh}. This is compatible with \cite[Ex.\,5.1]{CN0}, but cannot imply its proof immediately, since one has $\nu_{2d}^{(2)}\nes\nu_{2d}^{(3)}$ even in the case $f\eq z(x^4w{+}\1y^5{+}\1xy^4)$ (here the mapping cone of (\hl{3}{3}) seems to have two-dimensional support). Note also that in the case $\dim{\rm Sing}\,D\eq2$, one can prove (\hl{1.7.2}{1.7.2}) {\it for\1} $j\nes0,-1$ by an argument similar to the proof of Corollary~\hl{C1.4}{1.4}, since (\hl{1.4.12}{1.4.12}) can be extended for $\A_f^{n-2}$.
\bs\bs\htt{S2}{}
\vbox{\centerline{\bf 2. Proof of the main theorems}
\bsn
In this section, we prove the main theorems applying the assertions in the previous section.}
\par\htt{2.1}{}\msn
{\bf 2.1.~Proof of Theorem~\hl{T1}{1}.} By Proposition~\hl{P1.5}{1.5}, the stalk of the twisted logarithmic complex at $p$ is acyclic under the assumption (\hl{5}{5}). Indeed, the eigenvalues of the action of $L_{\xi_p}$ on the $\A_{f,p}^jf^{-1}$ are contained in $e_p^{-1}\1\Z$ with $e_p$ as in (\hl{4}{4}) by definition, whereas $\alt_p\notin e_p^{-1}\1\Z$ by the assumption (\hl{5}{5}). So it remains to show the vanishing
\htt{2.1.1}{}
$$H^j(V\stm D,L|_{V\setminus D})=0,
\leqno(2.1.1)$$
\par\nin where $V$ is a sufficiently small open neighborhood of $p$ in $X$. Using the $\C^*$-action, we may assume $V\eq\C^n$ with $D$ defined by a weighted homogeneous polynomial $f_p\eq\mprod_k\,f_{p,k}$ in $\C^n$.
\sk
Set
\vskip-7mm
$$a_i:=e_pw_{p,i}\q(i\in[1,n]).$$
\par\nin \vskip-.2mm\nin
Consider the finite map
$$\pi:V':=\C^n\ni(x_1,\dots,x_n)\mapsto(x_1^{a_1},\dots,x_n^{a_n})\in V\eq\C^n.$$
\par\nin Put
\vskip-8mm
$$g:=\pi^*\!f_p,\q g_k:=\pi^*\!f_{p,k}.$$
\par\nin Then $g,g_k$ are homogeneous polynomials of degree $d:=e_p$ and $d_k:=e_pd_{p,k}$ respectively.
\sk
Let $\rho:\Vt\to V'$ be the blow-up along $0\in V'$.
Let $D',L'$ be the inverse image or pullback of $D,L$ by $\pi$, and similarly for $\Dt,\Lt$ replacing $\pi$ with $\pit:=\pi\ssc\rho$. Here $\Lt$ is identified with $L'$.
Since $L$ is a direct factor of the direct image of $L'$ by the finite morphism $\pi$, it is sufficient to show that
\htt{2.1.2}{}
$$H^j(\Vt\stm\Dt,\Lt|_{\Vt\setminus\Dt})=0.
\leqno(2.1.2)$$
\par\nin \sk
Let $\la_0$ be the local monodromy (eigenvalue) of the local system $\Lt$ around the exceptional divisor $E$ of $\rho$. Since $\Vt\stm\Dt$ is a $\C^*$-bundle over $\PP^{n-1}\stm\{g\eq0\}$, the vanishing (\hl{2.1.2}{2.1.2}) is reduced to the inequality
\htt{2.1.3}{}
$$\la_0\ne 1,
\leqno(2.1.3)$$
\par\nin using the Leray-type spectral sequence.
\sk
In order to calculate $\la_0$, consider a generic line $\ell$ in $V'$ passing through the origin. Its pullback to $\Vt$ intersects $E$ transversally. We perturb this $\ell$ slightly so that it intersects $D'$ transversally at smooth points. The above $\la_0$ is the product of the local monodromy eigenvalues for all the intersection points. (Note that $\ell$ is a complex line.) Let $\la_k$ be the contribution coming from the intersection points with $\{g_k\}$ so that $\la_0\eq\mprod_k\,\la_k$. In view of the hypothesis (\hl{5}{5}), the inequality (\hl{2.1.3}{2.1.3}) is then reduced to the following.
\htt{2.1.4}{}
$$\la_k=e^{-2\pi i\1d_k\1\al_{p,k}}.
\leqno(2.1.4)$$
\par\nin \sk
In the case $f_{p,k}$ does not coincide with any $x_i$ (up to non-zero constant multiple), the intersection $\{g_k\eq0\}\caps(\C^*)^n$ is non-empty, and we get (\hl{2.1.4}{2.1.4}), since $\pi$ is unramified over $(\C^*)^n$ and the intersection number of $\{g_k\eq0\}$ and $\ell$ coincides with the degree of $g_k$, that is, $d_k$.
\sk
In the other case, we have $f_{p,k}\eq x_i$ (up to non-zero constant multiple) for some $i$. The intersection number of $\{x_i\eq0\}$ with $\ell$ is 1, and the local monodromy of $L'$ around $\{x_i\eq0\}$ is given by $\la_i^{a_i}$ with $\la_i\eq e^{-2\pi i\al_{p,k}}$ and $d_{p,k}\eq w_{p,i}$ (since $f_{p,k}\eq x_i$). So (\hl{2.1.4}{2.1.4}) follows. This finishes the proof of Theorem~\hl{T1}{1}.
\par\htt{2.2}{}\msn
{\bf 2.2.~Proof of Theorem~\hl{T2}{2}.} We have the commutative diagram
\htt{2.2.1}{}
$$\begin{array}{cccccccc}\Hc^n\bl(\DRl_X^{-1}(\M^{\ssb}_{X,\log}(L),\naal)\br)\,\,\longrightarrow\,\,\Hc^n\bl(\DRl_X^{-1}(\M^{\ssb}_X(L),\naal)\br)\\ |\!|\q\q\q\q\q\q\q\q\q\q\q\q\q\q\q |\!|\raise5mm\h{}\raise-2mm\h{}\\ \D_X/\D_X\Tht_{f,\al-1}\q\buildrel{\rho}\over\onto\q\D_Xf^{\al-1}\q\buildrel{\iota}\over\into\q\OO_X(*D)f^{\al-1}\end{array}
\leqno(2.2.1)$$
\par\nin see Remark~\hl{R1.7a}{1.7a}. If condition~(\hl{a2}{a}) holds, then the upper horizontal morphism is bijective, hence the lower horizontal morphisms $\rho,\iota$ are. Thus conditions (\hl{c2}{c}) and (\hl{f2}{f}) hold. We can verify the equivalence (\hl{c2}{c})\,$\Leftrightarrow$\,(\hl{d2}{d}) using the definitions of $\Tht_{f,\al}$ and $\Theta_X(-\log D)$, see Remark~\hl{R1.7a}{1.7a}. The equivalence (\hl{e2}{e})\,$\Leftrightarrow$\,(\hl{f2}{f}) follows for instance from \cite[Thm.\,1]{rp}, see Proposition~\hl{P2.2}{2.2} below and also \cite{Bu}.
\sk
We may then assume that the divisor $D$ is everywhere positively weighted homogeneous and the complex $\DRl^{-1}_X\bl(\M^{\ssb}_{X,\log}(L),\naal\br)[n]$ is quasi-isomorphic to a $\D_X$-module. The first assumption implies the equivalence (\hl{a2}{a})\,$\Leftrightarrow$\,(\hl{b2}{b}) by Lemma~\hl{L1.6}{1.6} and Proposition~\hl{P1.6}{1.6}. From the second one, we can deduce that
\htt{2.2.2}{}
$$\Hc^j\Cc^{\ssb}\eq0\q\h{if}\q j\ne n{-}1,n,
\leqno(2.2.2)$$
\par\nin with $\Cc^{\ssb}$ the mapping cone of the comparison morphism (\hl{2}{2}) applied by $\DRl_X^{-1}$. (This follows from the long exact sequence associated with the mapping cone.)
\sk
Assume condition~(\hl{c2}{c}) holds (that is, $\rho$ is injective), but (\hl{a2}{a}) does not. We then get that 
\htt{2.2.3}{}
$$\Hc^n\Cc^{\ssb}\ne0,\q\Hc^j\Cc^{\ssb}\eq0\,\,\,(\forall\,j\ne n).
\leqno(2.2.3)$$
\par\nin We may assume that the cohomology $\D$-module $\Hc^n\Cc^{\ssb}$ is supported at a point $p\ins D$ using Lemma~\hl{L1.6}{1.6} and cutting $D$ by a transversal slice to a maximal-dimensional stratum contained in the support. Then this $\D$-module is a finite direct sum of $\B_p\eq\C[\dd_{x_1},\dots,\dd_{x_n}]$ supported at $p$ (using a special case of Kashiwara's equivalence as in the proof of Proposition~\hl{P1.6}{1.6}). This implies that the mapping cone of the comparison morphism (\hl{2}{2}) is acyclic except at the highest degree $n$, since $\DR_X(\B_p)=\C_{\{p\}}$. However, this contradicts Corollary~\hl{C1.5}{1.5}. We thus get that (\hl{c2}{c})\,$\Rightarrow$\,(\hl{a2}{a}).
\sk
We can prove the implication (\hl{f2}{f})\,$\Rightarrow$\,(\hl{a2}{a}) by an argument similar to the above one, where $n$ is replaced by $n{-}1$ in (\hl{2.2.3}{2.2.3}). This finishes the proof of Theorem~\hl{T2}{2}.
\par\htt{P2.2}{}\msn
{\bf Proposition~2.2.} {\it For a holomorphic function $f$ on a complex manifold $X$, we have}
\htt{2.2.4}{}
$$\aligned&\D_Xf^{\al}=\D_Xf^{\al+1}\q\h{\it if}\q b_f(\al)\nes0,\\&\D_Xf^{\al}\ne\D_Xf^{\al+1}\q\h{\it if}\q b_f(\al)\eq0,\,\,b_f(\al{-}k)\ne0\,\,(\forall\,k\ins\Z_{>0}).\endaligned
\leqno(2.2.4)$$
\par\nin \msn
{\it Proof.} By \cite[Thm.\,1]{rp} there are regular holonomic $\D_X$-modules $\M_f^{\beta}$ endowed with a finite increasing filtration $G_{\ssb}$ and an nilpotent endomorphism $N$ for $\beta\in(0,1]$ such that $G_{\ssb}$ is stable by $N$ and
\htt{2.2.5}{}
$$\aligned\Gr^G_j\M_f^{\beta}\eq0&\iff b_f(-\beta{-}j)\nes0,\\ \Gr^G_j(\M_f^{\beta}/N\M_f^{\beta})\eq0&\iff\D_Xf^{-\beta-j}\eq\D_Xf^{-\beta-j+1}.\endaligned
\leqno(2.2.5)$$
\par\nin (Here $N$ is not {\it strictly compatible\1} with $G$ in general.) The first assertion then follows. One can also get this by setting $s\eq\al$ in the functional equation associated with $b_f(s)$.
\sk
For the second assertion, we see that the second hypothesis implies the equalities
\htt{2.2.6}{}
$$\D_Xf^{\al}\eq\D_Xf^{\al-k}\,\,(\forall\,k\ins\Z_{>0}),\q\M_f^{\beta}\eq G_j\M_f^{\beta},
\leqno(2.2.6)$$
\par\nin where $\al\eq{-}\beta{-}j$. From the second equality of (\hl{2.2.6}{2.2.6}) we can deduce the isomorphism
\htt{2.2.7}{}
$$\Gr^G_j(\M_f^{\beta}/N\M_f^{\beta})=\Gr^G_j\M_f^{\beta}/N\Gr^G_j\M_f^{\beta}.
\leqno(2.2.7)$$
\par\nin So the assertion follows (since the cokernel of a nilpotent endomorphism of a nonzero object of an abelian category is nonzero). This finishes the proof of Proposition~\hl{P2.2}{2.2}.
\par\htt{R2.2}{}\msn
{\bf Remark~2.2.} In the case the divisor $D$ is everywhere positively weighted homogeneous and the complex $\DRl_X^{-1}\bl(\M^{\ssb}_{X,\log}(L),\naal\br)[n]$ is quasi-isomorphic to a $\D_X$-module, this can be shown using Theorem~\hl{T1.3}{1.3}, Proposition~\hl{P1.5}{1.5}, and Remark~\hl{R1.2}{1.2} by induction on strata. Indeed, the argument is quite similar to the proof of Theorem~\hl{T3}{3} below. Here we have to replace the unipotent monodromy part in (\hl{2.3.1}{2.3.1}) below by the $\lambda$-eigenspace with $\lambda\eq e^{2\pi i\al}$, and the degree 1 part in (\hl{2.3.2}{2.3.2}) below by the degree $1{-}\al$ part. Note that $\om_p\wedge$ vanishes, see also Remark~\hl{R3.2e}{3.2e} below.
\par\htt{2.3}{}\msn
{\bf 2.3.~Proof of Theorem~\hl{T3}{3}.} We first show the assertion (ii). Under the hypothesis on the integral roots of $b_{f,q}(s)$ for $q\ne p$, the BS polynomial $b_f(s)$ has an integral root which is strictly smaller than $-1$ if and only if
\htt{2.3.1}{}
$$P^{n-1}H^{n-1}(F_{\!f},\C)_1\ne H^{n-1}(F_{\!f},\C)_1,
\leqno(2.3.1)$$
\par\nin see Theorem~\hl{T1.3}{1.3} (that is, \cite[Thm.\,2]{bCM}). By (\hl{1.2.4}{1.2.4}--\hl{1.2.5}{5}), the left-hand side of (\hl{2.3.1}{2.3.1}) is identified with the image of
\htt{2.3.2}{}
$$(\Hc_{f,p}^n)_1\to(\G_{f,p}^n)_1.
\leqno(2.3.2)$$
\par\nin Here we can replace $\G_{f,p}^n$ with $\G_{f,p}^n[t^{-1}]$ (the localization by $t$), since we take $\Gr_V^1$. (Indeed, the kernel and cokernel of the localization morphism are unions of subgroups annihilated by $t^k$ ($k\gg0$), hence $\Gr_V^1$ does not change under the localization by $t$.) The morphism (\hl{2.3.2}{2.3.2}) is then identified with
\htt{2.3.3}{}
$$H^n(\A_{f,p}^{\ssb}f^{-1},\ddd)_0\to\rlap{\raise-10pt\h{$\,\,\,\scriptstyle r$}}\indlim\,H^n(\A_{f,p}^{\ssb}f^{-r},\ddd)_0.
\leqno(2.3.3)$$
\par\nin where $_0$ denotes the degree 0 part. (Note that the differential $\ddd$ on $\A_{f,p}^{\ssb}$ commutes with multiplication by $f^k$ for $k\ins\N$.) Since the meromorphic de Rham complex is the inductive limit of $(\Om_X^{\ssb}(\log D)f^{-r},\ddd)$, the assertion~(ii) now follows from Proposition~\hl{P1.4}{1.4}.
\sk
The assertion~(i) also follows from Proposition~\hl{P1.4}{1.4} using the coincidence of ${\rm Ker}\,\iota_{f,p}^n$ with the torsion (see \cite[Thm.\,1]{BS} or \hl{1.2}{1.2}), since the morphism (\hl{2.3.2}{2.3.2}) is identified with (\hl{2.3.3}{2.3.3}). This finishes the proof of Theorem~\hl{T3}{3}.
\par\htt{2.4}{}\msn
{\bf 2.4.~Proofs of Corollaries~\hl{C2}{2} and \hl{C3}{3}.} Corollary~\hl{C2}{2} follows from Theorem~\hl{T3}{3} using \cite[Thm.\,2]{wh}, \cite[Cor,\,4.7]{DS2} as is explained after Corollary~\hl{C2}{2}, since the morphism (\hl{2.3.2}{2.3.2}) is identified with (\hl{2.3.3}{2.3.3}) in \hl{2.3}{2.3}.
\sk
For the proof of Corollary~\hl{C3}{3}, it is enough to show that $M_{kd}^{(2)}\ne 0$ for some $k\ins\Z_{\ges2}$ in the notation of Remark~\hl{R1.3a}{1.3a}, using Proposition~\hl{P1.4}{1.4}. This non-vanishing follows from {\it injectivity\1} of the morphism
\htt{2.4.1}{}
$$M'\into M^{(2)}\,\,\,\h{(or equivalently,}\,\,\,M'\cap{\rm Im}\,\ddd_1\eq0),
\leqno(2.4.1)$$
\par\nin (see \cite[Thm.\,3]{wh}) combined with a {\it symmetry\1} of the $\mu'_k$ with center $nd/2$ (see \cite[Cor.\,1]{DS2}), where $n\gess 3$. Indeed, the hypothesis implies that $\mu'_{dn'}>0$, since
\htt{2.4.2}{}
$$\mu'_k\pl\mu''_k\eq\nu_k\pl\gamma_k\q\q(k\ins\Z),
\leqno(2.4.2)$$
\par\nin (see \cite[(0.3)]{DS2}) and $\mu''_k\less\tau_Z\eq\mu_Z$. Here
\htt{2.4.3}{}
$$\aligned&\q\,\,\msum_k\gamma_kt^k:=\bl(\msum_{k=1}^{d-1}\1t^k\br)\raise1pt\h{${}^n$}\q\bl(\h{in particular,}\,\,\,\gamma_d\eq\tbinom{d-1}{n-1}\br),\\&M'\,{:=}\,H^0_{\mm}M,\,\,\,M''\,{:=}\,M/M',\,\,\,\mu'_k\,{:=}\,\dim M'_k,\,\,\,\mu''_k\,{:=}\,\dim M''_k,\endaligned
\leqno(2.4.3)$$
\par\nin with $\mm\sst\C[x_1,\dots,x_n]$ the maximal ideal, $M$ is as in Remark~\hl{R1.3a}{1.3a}, and we have $\gamma_d\less\gamma_{dn'}$, see for instance \cite[(4.11.1)]{wh}. This finishes the proofs of Corollaries~\hl{C2}{2} and \hl{C3}{3}.
\par\htt{R2.4}{}\msn
{\bf Remark~2.4.} The morphism (\hl{7}{7}) can be injective even if the pole order spectral sequence does not degenerate at $E_2$, since we need the non-degeneration exactly at the degree $d$ part for the non-injectivity as is shown by Theorem~\hl{T3}{3}\,(i). Set, for instance,
$$f\eq x^4\pl y^3z\pl z^3w\pl xyzw.$$
\par\nin The associated projective hypersurface has a unique singular point which has type $T_{3,4,8}$ (not $T_{3,3,4}$) with Milnor number 14 and Tjurina number 13 (according to calculations by a small computer program for non-degenerate functions and also by Singular \cite{Sing}). Using a small computer program based on the algorithm explained in \cite{nwh}, the numerical data of the pole order spectral sequence are given as follows:
$$\scalebox{0.8}{$\begin{array}{rrrrrrrrrrrrrrrrrrrrrrrrr}
k:&4 &5 &6 &7 &8 &9 &10 &11 &12\\
\gamma_k:&1 &4 &10 &16 &19 &16 &10 &4 &1\\
\mu_k:&1 &4 &10 &16 &19 &17 &13 &13 &13\\
\mu^{\scriptscriptstyle(2)}_k:&1 &4 &8 &8 &7 &4 &1 &1 &1\\
\mu^{\scriptscriptstyle(3)}_k:&1 &4 &7 &7 &6 &3 & & &\\
\nu_k:& & & & & &1 &3 &9 &12\\
\nu^{\scriptscriptstyle(2)}_k:& & & & & &1 &1 &1 &\\
\nu^{\scriptscriptstyle(3)}_k:& & & & & &1 &1 &1 &\end{array}$}$$
\par\nin Here $\mu^{(r)}_k\eq\dim M^{(r)}_k$ (and similarly for $\nu^{(r)}_k$), see Remark~\hl{R1.3b}{1.3b}. In this case we have $\mu^{(2)}_4\eq\mu^{(3)}_4$. This implies that we may have the injectivity of (\hl{7}{7}). (Some more calculation would be needed to see if the spectral sequence degenerates at $E_3$.)
\sk
On the other hand, let
$$f\eq x^5\pl y^4z\pl x^3y^2\pl w^5,$$
\par\nin so that $f\eq g\pl w^5$ with $g\ins\C[x,y,z]$ treated in \cite{hil}, \cite[Ex.\,5.3]{wh}, \cite[Ex.\,5.1]{nwh}, see also Remark~\hl{R1.7d}{1.7d}. We have
$$\scalebox{0.8}{$\begin{array}{rrrrrrrrrrrrrrrrrrrrrrrrr}
k:&4 &5 &6 &7 &8 &9 &10 &11 &12 &13 &14 &15 &16\\
\gamma_k:&1 &4 &10 &20 &31 &40 &44 &40 &31 &20 &10 &4 &1\\
\mu_k:&1 &4 &10 &20 &31 &40 &45 &46 &45 &44 &44 &44 &44\\
\mu^{\scriptscriptstyle(2)}_k:&1 &3 &4 &6 &7 &7 &7 &6 &5 &4 &4 &4 &4\\
\mu^{\scriptscriptstyle(3)}_k:& &1 &1 &2 &3 &3 &3 &2 &1 & & & &\\
\nu_k:& & & & & & &1 &6 &14 &24 &34 &40 &43\\
\nu^{\scriptscriptstyle(2)}_k:& & & & & & & & & & &1 &2 &3\\
\nu^{\scriptscriptstyle(3)}_k:& & & & & & & & & & & & &\end{array}$}$$
\par\nin Here $\mu^{(2)}_5\ne \mu^{(3)}_5$. This implies the non-injectivity of (\hl{7}{7}) for $j\eq4$. (This non-injectivity holds also for $g$ with $j\eq3$, see Remark~\hl{R1.7d}{1.7d}.) The above calculation implies that the weighted degree 1 part of the Brieskorn module has nonzero torsion, so the logarithmic comparison theorem fails for the non-weighted-homogeneous isolated singularity defined by $x^5\pl y^4\pl x^3y^2\pl w^5\eq0$ by an argument similar to Remark~\hl{R1.7d}{1.7d}.
\par\htt{2.5}{}\msn
{\bf 2.5.~Proof of Proposition~\hl{P1}{1}.} The equivalence (\hl{a}{a})\,$\Leftrightarrow$\,(\hl{b}{b}) follows from Corollary~\hl{C2}{2} and Proposition~\hl{P1.4}{1.4} using \cite[Thm.\,0.1]{DS1} (see Remark~\hl{R1.3c}{1.3c}) and the inclusion $F\sst P$ on $H^1(F_{\!f,0},\C)$. Here we may assume that $p\eq 0\,({\in}\,D)$, since the assertion for $p\ins D\stm\{0\}$ is an easy consequence of Proposition~\hl{P1.4}{1.4}. (Note that 1 is the unique integral spectral number of reduced plane curve singularities.) We have $\nu_{kd}^{(2)}\eq 0$ for $k\gess 3$ using the inclusion $F\sst P$ together with the $E_2$-degeneration of spectral sequence (see \cite[Thm.\,2]{wh}), where the differential $\ddd$ shifts the degree by $-d$, since $\ddd f\wedge$ preserves it. This vanishing together with that of $\mu_{kd}^{(2)}$ for $k\gess 2$ (which is equivalent to condition~(\hl{b}{b}) at 0) implies surjectivity of (\hl{7}{7}) at 0 for any $j\ins\Z$ using Proposition~\hl{P1.4}{1.4}.
\sk
The equivalence (\hl{b}{b})\,$\Leftrightarrow$\,(\hl{c}{c}) follows from Corollary~\hl{C2}{2}. For (\hl{c}{c})\,$\Leftrightarrow$\,(\hl{d}{d}), we have symmetries of the $\mu'_k$ and $\delta''_k:=\mu''_k\mi\nu_{k+d}$ with center $3d/2$ and $d$ respectively (see \cite[(16)]{wh}). Note that, if the differential $\ddd_1\,{:}\,N_{k+d}\,{\to}\,M''_k$ is not bijective for some $k\gess 2d$ (where $\delta''_k\eq0$ using a symmetry), then we get that $\nu_{k+d}^{(2)}\ne 0$ by (\hl{2.4.1}{2.4.1}), but this contradicts the $E_2$-degeneration of spectral sequence using the inclusion $F\sst P$. The argument is similar for (\hl{c}{c})\,$\Leftrightarrow$\,(\hl{e}{e}) employing also \cite[Thm.\,4.1]{DiPo}, which says that $\mu'_{k-1}\less\mu'_k$ for $k\less 3d/2$ and $\mu'_{k+1}\less\mu'_k$ for $k\gess 3d/2$. This finishes the proof of Proposition~\hl{P1}{1}.
\par\htt{R2.5}{}\msn
{\bf Remark~2.5.} Assume $f\eq x^d\pl g(y,z)w$ with $g(y,z)$ a reduced homogeneous polynomial of degree $d{-}1$ in $y,z$ and $d\gess 3$. The unipotent monodromy part of the vanishing cycle complex $\varphi_{f,1}\Q_X$ vanishes (using for instance \cite[Rem.\,3.2]{low} and the Thom-Sebastiani type theorem \cite[Thm.\,0.8]{mic}). In the case $f\eq x^5\pl y^4w\pl z^4w$, the pole order spectral sequence can be calculated as below (using a small computer program explained in Remark~\hl{R2.4}{2.4}):
$$\scalebox{0.8}{$\begin{array}{rrrrrrrrrrrrrrrrrrrrrrrrr}
k:&4 &5 &6 &7 &8 &9 &10 &11 &12 &13 &14 &15 &16\\
\gamma_k:&1 &4 &10 &20 &31 &40 &44 &40 &31 &20 &10 &4 &1\\
\mu_k:&1 &4 &10 &20 &31 &41 &48 &51 &52 &52 &52 &52 &52\\
\mu^{\scriptscriptstyle(2)}_k:& & &3 &3 &3 &3\\
\nu_k:& & & & & &1 &4 &11 &21 &32 &42 &48 &51\\
\nu^{\scriptscriptstyle(2)}_k:& & & & & & & &4 &4 &4 &4\\
\end{array}$}$$
\par\nin \bs\bs\htt{S3}{}
\vbox{\centerline{\bf 3. Hyperplane arrangement case}
\bsn
In this section, we give a simple proof of a stronger version of the comparison theorem for hyperplane arrangements as an immediate corollary of \cite{Sc}, \cite{DeSi}.}
\par\htt{3.1}{}\msn
{\bf 3.1.~Castelnuovo-Mumford regularity.} For a finitely generated graded $R$-module $M$ with $R:=\C[x_1,\dots,x_n]$, the {\it Castelnuovo-Mumford regularity\1} can be defined as
\htt{3.1.1}{}
$${\rm reg}\,M:=\max_{j,k}\,(c_{j,k}\mi j),
\leqno(3.1.1)$$
\par\nin by taking a {\it minimal graded free resolution}
\htt{3.1.2}{}
$$\to F_j\to\cdots\to F_0\to M\to 0,
\leqno(3.1.2)$$
\par\nin with $F_j=\mopl_k\,R(-c_{j,k})$ ($c_{j,k}\ins\Z,\,j\ins\N$), see for instance \cite[\S 4A]{Ei}.
\sk
By H.~Derksen and J.~Sidman, we have the following.
\par\htt{P3.1}{}\msn\vbox{\nin
{\bf Proposition~3.1} (\cite[Cor.\,3.7]{DeSi}). {\it Let $F$ be a free graded $R$-module which is freely generated in degree $m\,\,($that is, $F$ is isomorphic to a direct sum of copies of $R(-m))$. Let $M,M_i\sst F$ be graded $R$-submodules $(i\ins[1,n])$. Assume
\htt{3.1.3}{}
$$x_i\1M_i\subset M\subset M_i\q(i\ins[1,n]),
\leqno(3.1.3)$$
\par\nin and ${\rm reg}\,M_i\less r\mi 1\,\,(i\ins[1,n])$ for some $r\gess m\pl 2$. Then ${\rm reg}\,M\less r$.}}
\ms
By the same argument as in the proof of \cite[Prop.\,1.3]{4var}, which was inspired by \cite{Sc} and a remark before \cite[Thm.\,5.5]{DeSi} (where the degree is shifted by 1), we can immediately get the following (compare to \cite{Ba}, where the proof is formulated in a rather complicated manner).
\par\htt{C3.1}{}\msn\vbox{\nin
{\bf Corollary~3.1.} {\it Let $D$ be an essential reduced central hyperplane arrangement in $X:=\C^n$. Then we have the estimate
\htt{3.1.4}{}
$${\rm reg}\,\Ga\bl(X,\Om_X^j(\log D)\br)\les 0\q\q(\forall\,j\ins\Z),
\leqno(3.1.4)$$
\par\nin where algebraic logarithmic forms are used.}}
\msn
{\it Proof.} We argue by induction on $n$ and $d:=\deg D$ as in the proof of \cite[Prop.\,1.3]{4var}. Here we may assume $d>n$ (since the normal crossing case with $d\eq n$ is trivial) and $j\ins[1,n{-}1]$ (since $\Om_X^j(\log D)\eq\OO_X$ or $\OO_X(d{-}n)$ if $j\eq0$ or $n$). We may assume also that the coordinate hyperplanes $D_i:=\{x_i\eq0\}$ are contained in $D$ changing the coordinates if necessary (since $D$ is assumed to be essential). Let $D^{(i)}$ be the closure of $D\stm D_i$ in $X$. There are inclusions
\htt{3.1.5}{}
$$\Om_X^j(\log D^{(i)})\subset \Om_X^j(\log D)\subset x_i^{-1}\Om_X^j(\log D^{(i)})\q(j\ins\Z,\,i\in[1,n]).
\leqno(3.1.5)$$
\par\nin (It is enough to verify these at {\it smooth\1} points of $D$ using the Hartogs-type theorem for logarithmic forms.) If the $D^{(i)}$ are all essential, then the assertion follows from the inductive hypothesis applying Proposition~\hl{P3.1}{3.1} to
\htt{3.1.6}{}
$$M_i:=x_i^{-1}\Ga\bl(X,\Om_X^j(\log D^{(i)})\br)\q(i\ins[1,n]),
\leqno(3.1.6)$$
\par\nin where $F:=\Ga\bl(X,\Om_X^j(D)\br)$ with $m\eq j\mi d\les -2$ (since $j<n<d$).
\sk
Assume $D^{(n)}$ is not essential (replacing the order of coordinates if necessary). Then we have $f=x_n\1g$ with $g\in R':=\C[x_1,\dots,x_{n-1}]$ (changing the coordinates if necessary). Set
$$D':=g^{-1}(0)\subset X':=\C^{n-1},\q D'':=\{0\}\subset X'':=\C.$$ We have the K\"unneth formula (see for instance \cite{CDFV})
\htt{3.1.7}{}
$$\Om_X^j(\log D)=\Om_{X'}^j(\log D')\boxtimes\OO_{X''}\oplus
\Om_{X'}^{j-1}(\log D')\boxtimes\Om_{X''}^1(D'').
\leqno(3.1.7)$$
\par\nin (This also follows from the Hartogs-type theorem.) Here $\Om_{X''}^1(D'')$ is identified with a free graded $\C[x_n]$-module generated freely by $\ddd x_n/x_n$ (which has degree 0) taking the global sections.
So the assertion follows from the inductive hypothesis. We can thus proceed by induction on $n,d$. This finishes the proof of Corollary~\hl{C3.1}{3.1}.
\par\htt{R3.1}{}\msn
{\bf Remark~3.1.} The estimate of regularity (\hl{3.1.4}{3.1.4}) for $j\eq n{-}1$ is equivalent to the one for logarithmic vector fields \cite[Prop.\,1.3]{4var}. Indeed, there is an isomorphism induced by the interior product
\htt{3.1.8}{}
$$\Theta_X(-\log D)\ni\xi\,\mapsto\,f^{-1}\iota_{\xi}\1\om_0\ins\Om_X^{n-1}(\log D),
\leqno(3.1.8)$$
\par\nin with $\om_0:=\ddd x_1{\wedge}\cdots{\wedge}\ddd x_n$. Here the degree is shifted by $n\mi d$.
\par\htt{3.2}{}\msn
{\bf 3.2.~Higher cohomology vanishing.} From Corollary~\hl{C3.1}{3.1}, we can deduce the following.
\par\htt{T3.2}{}\msn
{\bf Theorem~3.2.} {\it Let $D$ be an essential reduced central hyperplane arrangement in $X:=\C^n$. Let $\pi:\Xt\to X$ be the blow-up at $0\ins X$ with $E:=\pi^{-1}(0)\,(\cong\PP^{n-1})$ the exceptional divisor. Setting $\Dt:=\pi^{-1}(D)$, we have the higher cohomology vanishing}
\htt{3.2.1}{}
$$H^i\bl(E,\Om_{\Xt}^j(\log\Dt){\otimes}_{\OO_{\Xt}}\OO_E(k)\br)=0\q(i\,{>}\,0,\,k\gess{-}1,\,j\ins\Z).
\leqno(3.2.1)$$
\par\nin \msn
{\it Proof.} Set
$$M^j:=\Ga\bl(X,\Om_X^j(\log D)\br)\q(j\ins\Z).$$
\par\nin The graded $R$-module $M^j$ corresponds to the $\OO_X$-module $\Om_X^j(\log D)$, and also to the $\OO_E$-module $\F^j$ such that
\htt{3.2.2}{}
$$\Ga\bl(E,\F^j(k)\br)=M^j_k\q(\forall\,k\ins\Z).
\leqno(3.2.2)$$
\par\nin Indeed, we have the isomorphisms at least for $k\gg 0$ (since $E\eq{\rm Proj}\,R$, see \cite{Ha}). The graded $R$-module $\mopl_{k\in\Z}\,\Ga\bl(E,\F^j(k)\br)$ corresponds to the direct image by the inclusion $X\stm\{0\}\into X$, and we can apply the Hartogs-type extension theorem.
\sk
We have moreover the isomorphisms
\htt{3.2.3}{}
$$\F^j(k)=\Om_{\Xt}^j(\log\Dt){\otimes}_{\OO_{\Xt}}\,\I_E^k/\I_E^{k+1}\q(k\gess 0),
\leqno(3.2.3)$$
\par\nin where $\I_E\sst\OO_{\Xt}$ is the ideal sheaf of $E\sst\Xt$. Indeed, the right-hand side of (\hl{3.2.3}{3.2.3}) can be identified with the subsheaf of $\Om_{\Xt}^j(\log\Dt)|_E$ on which the action of the Lie derivation $L_{\xit}$ is given by multiplication by $k$ with $\xit$ the {\it pull-back\1} of the vector field
$$\xi:=\msum_{i=1}^n\,x_i\dd_{x_i},$$
\par\nin under the {\it birational\1} morphism $\pi$. Here we use locally the analytic K\"unneth formula for logarithmic forms as in (\hl{3.1.7}{3.1.7}) around each point of $E\sst\Xt$ together with GAGA. Note that the action of $L_{\xi}$ on $M^j_k$ is given by multiplication by $k$.
\sk
Since the normal bundle of $E\sst\Xt$ is $\OO_E(-1)$, we get that
\htt{3.2.4}{}
$$\I_E^k/\I_E^{k+1}\cong\OO_E(k)\q(k\gess 0).
\leqno(3.2.4)$$
\par\nin The vanishing (\hl{3.2.1}{3.2.1}) now follows from Corollary~\hl{C3.1}{3.1} taking a minimal graded free resolution as in (\hl{3.1.2}{3.1.2}), since we have for $i\,{>}\,0$
\htt{3.2.5}{}
$$H^i\bl(E,\OO_E(k)\br)\eq 0\q\h{unless}\q i\eq n{-}1,\,k\les-n.
\leqno(3.2.5)$$
\par\nin This finishes the proof of Theorem~\hl{T3.2}{3.2}.
\ms
Set $\al_Z\,{:=}\,\msum_{D_k\supset Z}\,\al_k$ for $Z\sst X$ in the notation of the introduction. Let $\M^{\ssb}_{X,\log}(L)_{\ges m}$ be the subcomplex whose quotient is supported at $0$ and whose stalk at $0$ is generated by local sections which are eigenvectors of the Lie derivation $ L_{\xi}$ of the Euler field $\msum_i\,x_i\dd_{x_i}$ with eigenvalues at least $m$. 
As a corollary of Theorem~\hl{T3.2}{3.2}, we get the following.
\par\htt{C3.2a}{}\msn\vbox{\nin
{\bf Corollary~3.2a.} {\it In the notation of Theorem~{\rm\hl{T3.2}{3.2}}, let $m\ins\Z_{\ges-1}$. Set $X^*\,{:=}\,X\stm\{0\}$ with $j_0\,{:}\,X^*\,{\into}\,X$ the inclusion. We have the isomorphism in $D^b_c(X,\C)\,{:}$
\htt{3.2.6}{}
$$\M^{\ssb}_{X,\log}(L)_{\ges m}=\begin{cases}\R(j_0)_*\M^{\ssb}_{X^*,\log}(L)&\h{if}\,\,\,\al_{\{0\}}{+}m\notin\Z_{\ges 1},\\(j_0)_!\M^{\ssb}_{X^*,\log}(L)&\h{if}\,\,\,\al_{\{0\}}{+}m\notin\Z_{\les 0},\end{cases}
\leqno(3.2.6)$$
\par\nin where $\naal$ is omitted to simplify the notation}}
\msn
{\it Proof.} Let $\pi:\Xt\to X$, and $\Dt$, $E$ be as in Theorem~\hl{T3.2}{3.2}. Set
$$\Lc_{\Xt}:=\pi^*\Lc_X(-mE).$$
\par\nin Here $\Lc_X$ is trivialized by using the twisted differential as in the introduction, and $(-mE)$ means the tensor product with $\OO_{\Xt}(-mE)$. The residue $\al_E$ along the exceptional divisor $E$ is given by
\htt{3.2.7}{}
$$\al_E=\al_{\{0\}}\pl m\,\,\bl(=\msum_{k=1}^d\,\al_k\pl m\br).
\leqno(3.2.7)$$
\par\nin This can be shown by taking the pull-backs of $\ddd f_k/\!f_k$ to $\Xt$, where the $f_k$ are defining linear functions of $D_k$. (Note that $+m$ on the right-hand side of (\hl{3.2.7}{3.2.7}) corresponds to $(-mE)$ in the definition of $\Lc_{\Xt}$.) Define $\M_{\Xt,\log}^{\ssb}(L)$ using this $\Lc_{\Xt}$. Let $\jt:\Xt\stm E\into\Xt$ be the inclusion. We have the isomorphisms in $D^b_c(\Xt,\C)$\,:
\htt{3.2.8}{}
$$\M_{\Xt,\log}^{\ssb}(L)=\begin{cases}\R\jt_*\M^{\ssb}_{X^*,\log}(L)&\h{if}\,\,\,\al_{\{0\}}{+}m\notin\Z_{\ges 1},\\ \jt_!\M^{\ssb}_{X^*,\log}(L)&\h{if}\,\,\,\al_{\{0\}}{+}m\notin\Z_{\les 0},\end{cases}
\leqno(3.2.8)$$
\par\nin using the K\"unneth formula, see (\hl{3.1.7}{3.1.7}). Indeed, there is a well-known isomorphism
\htt{3.2.9}{}
$$C(t\dd_t\pl\beta\,{:}\,\OO_{\Delta}\,{\to}\,\OO_{\Delta})=\begin{cases}\R j'_*L'[1]&\h{if}\,\,\,\beta\notin\Z_{\ges 1},\\ j'_!L'[1]&\h{if}\,\,\,\beta\notin\Z_{\les 0}.\end{cases}
\leqno(3.2.9)$$
\par\nin Here $\Delta$ is an open disk with coordinate $t$, and $L'$ is a local system on $\Delta^*:=\Delta\stm\{0\}$ with $j':\Delta^*\into\Delta$ the inclusion.
\sk
By Theorem~\hl{T3.2}{3.2} together with \cite[III, Thm.\,11.1]{Ha} (using (\hl{3.2.4}{3.2.4}) and GAGA), we get the vanishing
\htt{3.2.10}{}
$$R^i\pi_*\M_{\Xt,\log}^j(L)=0\q(i>0,\,j\ins\Z),
\leqno(3.2.10)$$
\par\nin since $\OO_{\Xt}(E){\otimes}_{\OO_{\Xt}}\OO_E\cong\OO_E(-1)$ and the range of $k$ in (\hl{3.2.1}{3.2.1}) is given by $k\gess{-}1$ (and we assume $m\ins\Z_{\ges-1}$). By (\hl{3.2.8}{3.2.8}) this implies the isomorphisms in $D^b_c(X,\C)$\,:
\htt{3.2.11}{}
$$\aligned&\pi_*\M_{\Xt,\log}^{\ssb}(L)=\R\pi_*\M_{\Xt,\log}^{\ssb}(L)\\
&\q=\begin{cases}\R(j_0)_*\M^{\ssb}_{X^*,\log}(L)&\h{if}\,\,\,\al_{\{0\}}{+}m\notin\Z_{\ges 1},\\(j_0)_!\M^{\ssb}_{X^*,\log}(L)&\h{if}\,\,\,\al_{\{0\}}{+}m\notin\Z_{\les 0}.\end{cases}\endaligned
\leqno(3.2.11)$$
\par\nin \sk
On the other hand, by the definition of $\Lc_{\Xt}$ just before (\hl{3.2.7}{3.2.7}), we can get the isomorphism of complexes
\htt{3.2.12}{}
$$\pi_*\M_{\Xt,\log}^{\ssb}(L)=\M_{X,\log}^{\ssb}(L)_{\ges m},
\leqno(3.2.12)$$
\par\nin using the actions of Lie derivations $L_{\pi^*\xi}$ and $L_{\xi}$ (together with the Hartogs-type extension theorem for logarithmic forms), where $\xi$ is as in the proof of Theorem~\hl{T3.2}{3.2}. So the isomorphism (\hl{3.2.6}{3.2.6}) follows. This completes the proof of Corollary~\hl{C3.2a}{3.2a}.
\sk
Corollary~\hl{C3.2a}{3.2a} implies a stronger version of the comparison theorem as follows (compare to \cite{Ba} where the argument is more complicated using local cohomology).
\par\htt{C3.2b}{}\msn
{\bf Corollary~3.2b.} {\it In the notation of Corollary~{\rm\hl{C3.2a}{3.2a}} assume $\al_k\notin\Z_{\ges 1}$ for any $k$, and $\al_Z\notin\Z_{\ges 2}$ for any dense edge $Z$. Then the comparison morphism~{\rm(\hl{2}{2})} is a $D$-quasi-isomorphism.}
\msn
{\it Proof.} We argue by induction on the codimension of strata of the stratification associated to the hyperplane arrangement $D$. The assertion is well known in the codimension~1 case, see \cite{De}. Assume the assertion is proved outside a given stratum locally. Here we may assume that the closure of the stratum is a dense edge. In the other case, we can apply the compatibility of the canonical morphism (\hl{3}{3}) with external product $\boxtimes$ using the K\"unneth formula and the inductive hypothesis. (The argument becomes more complicated if we consider only the morphism (\hl{2}{2}), since there is some difference between external products of $\OO$-modules and $\C$-complexes.)
\sk
Cutting $D$ by a transversal space passing through a general point of the stratum, we may assume that the stratum is $0$-dimensional. So $D$ is an essential indecomposable central hyperplane arrangement, and the assertion holds outside the origin.
\sk
Since $\al_{\{0\}}\notin\Z_{\ges 2}$, we get the quasi-isomorphism
\htt{3.2.13}{}
$$\M_{X,\log}^{\ssb}(L)_{\ges-1}\simto\M_{X,\log}^{\ssb}(L),
\leqno(3.2.13)$$
\par\nin using Proposition~\hl{P1.5}{1.5} (where the grading of $\A_{f,0}^{\ssb}$ is indexed by $\tfrac{1}{d}\1\Z$ with $d\,{:=}\,\deg f$ instead of $\Z$ and $\alt_0\eq\tfrac{1}{d}\,\al_{\{0\}}$). The isomorphism (\hl{3}{3}) then follows from Corollary~\hl{C3.2a}{3.2a} for $m\eq{-}1$ using Lemma~\hl{L1.6}{1.6}, Proposition~\hl{P1.6}{1.6} and the inductive hypothesis. This terminates the proof of Corollary~\hl{C3.2b}{3.2b}.
\par\htt{R3.2a}{}\msn
{\bf Remark~3.2a.} It is easy to see that the condition $\al_Z\notin\Z_{\ges2}$ is optimal in the case ${\rm codim}\,Z\eq 2$, see Example~\hl{E3.2}{3.2} below.
\sk
We can prove also the following (see \cite{CN} for the free divisor case).
\par\htt{C3.2c}{}\msn\vbox{\nin
{\bf Corollary~3.2c.} {\it Let $D\sst X$ be a reduced affine or projective hyperplane arrangement. For a dense edge $Z\sst D$, let $Z^o\sst Z$ be the complement of the union of hyperplanes of $D$ not containing $Z$. Set $\delta_Z:=\max_{j\in\Z}\delta_Z^{(j)}$ with
\htt{3.2.14}{}
$$\delta_Z^{(j)}:={\rm mult}_ZD-\min\bl\{k\ins\Z\mid(\A^j_{f,p})_k\ne0,\,p\ins Z^o\br\}\q(j\ins\Z).
\leqno(3.2.14)$$
\par\nin In the notation of Corollary~{\rm\hl{C3.2b}{3.2b},} assume $\al_Z\notin\Z_{\les\delta_Z}$ for any dense edge $Z\sst D$ $($where $\al_Z\eq\al_k$ with $\delta_Z\eq0$ if $Z\eq D_k)$. Then we have the canonical isomorphism in $D^b_{\rm hol}(\D_X)\,{:}$
\htt{3.2.15}{}
$$\DRl^{-1}_X\bl(\M^{\ssb}_{X,\log}(L),\naal\br)[n]\simto\DD\bl(\M_X(L^{\vee})\br).
\leqno(3.2.15)$$
\par\nin Here $L^{\vee}$ is the dual local system of $L$ and $\DD$ is the dual functor. In particular, there is a natural quasi-isomorphism
\htt{3.2.16}{}
$$(j_U)_!L\simto\bl(\M^{\ssb}_{X,\log}(L),\naal\br)\q\h{in}\,\,\,D^b_c(X,\C),
\leqno(3.2.16)$$
\par\nin or equivalently, $\Hc^j\bl(\M^{\ssb}_{X,\log}(L),\naal\br){}_p\eq0$ for any $p\ins D$, $j\ins\Z$.}}
\msn
{\it Proof.} In the definition (\hl{3.2.14}{3.2.14}), the multiplicity ${\rm mult}_ZD$ of $D$ along $Z$ coincides with the number of hyperplanes in $D$ containing $Z$, and the grading of $\A^j_{f,p}$ is indexed by $\Z$ (counted by using the vector field $\xi$ as in the proof of Theorem~{\rm\hl{T3.2}{3.2}}). We have
\htt{3.2.17}{}
$$\delta_Z\ins\Z_{\ges1}\q\h{unless}\q{\rm codim}_XZ\eq1.
\leqno(3.2.17)$$
\par\nin This is verified by reducing to the {\it indecomposable essential\1} central arrangement case, and considering the highest forms (that is, $j\eq n)$.
\sk
The proof of Corollary~\hl{C3.2c}{3.2c} is quite similar to that of Corollary~\hl{C3.2b}{3.2b}. We argue by induction on the codimension of strata. In the codimension 1 case, we have the isomorphism noted after (\hl{3.2.8}{3.2.8}).
\sk
For the inductive argument, we may assume that $D$ is an essential central arrangement with $Z\eq\{0\}$. Since $\al_{\{0\}}\notin\Z_{\les\delta_{\{0\}}}$, we can prove the quasi-isomorphism
\htt{3.2.18}{}
$$\M_{X,\log}^{\ssb}(L)_{\ges0}\simto\M_{X,\log}^{\ssb}(L),
\leqno(3.2.18)$$
\par\nin using Proposition~\hl{P1.5}{1.5}. The assertion then follows from Corollary~\hl{C3.2a}{3.2a} for $m\eq0$ using an argument similar to the proof of Proposition~\hl{P1.6}{1.6} and the inductive hypothesis. This completes the proof of Corollary~\hl{C3.2c}{3.2c}.
\par\htt{R3.2b}{}\msn
{\bf Remark~3.2b.} Corollary~\hl{C3.2c}{3.2c} is {\it optimal\1} in the following sense: Assume there is a dense edge $Z\sst D$ with $\al_Z\eq\delta_Z$ and $\al_{Z'}\notin\Z_{\les\delta_{Z'}}$ for any dense edge $Z'\sst D$ strictly containing $Z$. Then the morphism (\hl{3.2.16}{3.2.16}) is not a quasi-isomorphism at any $p\in Z^o$, since (\hl{3.2.18}{3.2.18}) is not. This argument may become quite complicated if the condition $\al_Z\eq\delta_Z$ is replaced by $\al_Z\ins\Z_{\les\delta_Z}$ and we have $\al_Z\less\delta^{(j)}_Z$ for several $j\ins\Z$, see Example~\hl{E3.2}{3.2} below for a simple case.
\par\htt{R3.2c}{}\msn
{\bf Remark~3.2c.} It is easy to verify that
\htt{3.2.19}{}
$$\delta_Z\less{\rm mult}_ZD\mi 3,
\leqno(3.2.19)$$
\par\nin by reducing to the case of an essential central arrangement with $Z\eq\{0\}$. Here $\delta_{\{0\}}^{(j)}\less d\mi j$ ($j\ins[2,\dim X]$) with $d\,{:=}\,\deg D\eq\deg f$ by definition (\hl{3.2.14}{3.2.14}), and we have $\delta_{\{0\}}^{(1)}\eq0$ since $D$ is reduced so that $\A_f^1=\OO_X\ddd f$, see for instance \cite[Prop.\,2]{JKSY2}. If $\delta_{\{0\}}\eq\delta_{\{0\}}^{(2)}\eq d\mi 2$, we look at $\bl(\A_{f,0}^2\br){}^{\h{}}_2$ after restricting $D$ to a sufficiently general 3-dimensional affine subspace containing $0$. (Note that $D$ is non-essential if there is a non-trivial linear relation with $\C$-coefficients among the $\dd_{x_i}f$.)
\par\htt{E3.2}{}\msn
{\bf Example~3.2.} Let $X\eq\C^2$ with $f$ a reduced homogeneous polynomial in two variables of degree $d\gess 3$. Assume $\al_{\{0\}}\eq\sum_{i=k}^d\al_k\ins\Z\cap[1,d{-}1]$. We have $\delta_{\{0\}}\eq d{-}2\eq{-}\chi(\PP^1\stm Z)$, and
\htt{3.2.20}{}
$$\dim H^j(X\stm D,L)=\begin{cases}d{-}2&\h{if}\,\,\,j\eq 1\,\,\h{or}\,\,2,\\ \1 0&\h{otherwise},\end{cases}
\leqno(3.2.20)$$
\par\nin using the Leray-type spectral sequence associated with a $\C^*$-bundle (see Remark~\hl{R1.4c}{1.4c}) or the Thom-Gysin sequence (see for instance \cite[Sect.\,1.3]{RSW}). Indeed, the local system $L$ is the pullback of a local system $L'$ of rank 1 on $\PP^1\stm Z$, since $\al_{\{0\}}\ins\Z$. The spectral sequence degenerates at $E_2$, since the {\it Euler class\1} vanishes by restricting it to $\PP^1\stm\{p\}$ for $p\ins Z$. Moreover the restriction of $\bl(\Om_X^{\ssb}(\log D),\naal\br)$ to the complement of $0\ins X$ is quasi-isomorphic to the pullback of $\bl(\Om_{\PP^1}^{\ssb}(\log Z),\naal{}'\br)$, which is quasi-isomorphic to $K':=\R j''_*j'_!L'$. Here $j\eq j''\ssc j':\PP^1\stm Z\into\PP^1$ is a factorization depending on the $\al_k$ ($k\ins[1,d]$), and we may assume $j',j''$ are not the identity, since $\al_{\{0\}}\eq\sum_{k=1}^d\al_k\ins[1,d{-}1]$. We then get that $H^i(\PP^1,K')\eq0$ ($i\nes1$), and $\dim H^1(\PP^1,K')\eq d{-}2$, since $\chi(K')\eq\chi(\R j_*L')\eq 2{-}d$.
\sk
On the other hand, Proposition~\hl{P1.5}{1.5} implies that
\htt{3.2.21}{}
$$\dim\Hc^j\bl(\Om_X^{\ssb}(\log D),\naal\br){}_0=\begin{cases}d{-}1{-}\al_{\{0\}}&\h{if}\,\,\,j\eq 1\,\,\h{or}\,\,2,\\ \1 0&\h{otherwise}.\end{cases}
\leqno(3.2.21)$$
\par\nin Indeed, $\A_f^1\eq\OO_X\ddd f$ (since $D$ has an isolated singularity), and hence $(\A_{f,0}^1)_k\eq0$ for $k\,{<}\,d$.
\sk
So the morphisms (\hl{2}{2}) and (\hl{3.2.16}{3.2.16}) are both {\it non-quasi-isomorphisms\1} if $\al_{\{0\}}\ins\Z\cap[2,d{-}2]$. Note that the {\it higher vanishing\1} is irrelevant to the {\it twist\1} of the differential operator (unless one tries to use the direct image of the {\it canonical\1} Deligne extension whose residues are contained either in $[0,1)$ or in $(0,1]$). If one applies Corollary~\hl{C3.2c}{3.2c} with $m\eq0$ assuming $\al_{\{0\}}\ins\Z_{>0}$ with no condition on each $\al_k$, then one gets the zero-extension by the inclusion $X\stm\{0\}\into X$, but (\hl{3.2.18}{3.2.18}) does not necessarily hold unless $\al_{\{0\}}\gess d{-}1$.
\sk
Note finally that Corollary~\hl{C3.2b}{3.2b} implies the following generalization of \cite[Thm.\,5.3]{Wa2} (where $\al\eq0$).
\par\htt{C3.2d}{}\msn
{\bf Corollary~3.2d.} {\it Let $\al\ins\C\stm\bl(\Z_{\ges 1}\cup\mcup_Z\,(m_Z)^{-1}\Z_{\ges 2}\br)$. Here $Z$ runs over the dense edges of $D$, and $m_Z$ is the number of hyperplanes of $D$ containing $Z$. Then the annihilator ${\rm Ann}_{\D_X}(f^{\al-1})$ of $f^{\al-1}$ in $\D_X$ is generated by $\Tht_{f,\al-1}$ in the notation of Remark~{\rm\hl{R1.7a}{1.7a}}.}
\msn
{\it Proof.} Set $\al_k:=\al$ ($\forall\,k$). The hypothesis of Corollary~\hl{C3.2b}{3.2b} is satisfied by the assumption on $\al$. So the assertion follows from this corollary together with Theorem~\hl{T2}{2} in the general case.
\par\htt{R3.2d}{}\msn
{\bf Remark~3.2d.} The hypothesis of Corollary~\hl{C3.2d}{3.2d} implies that $b_f(\al{-}1{-}j)\nes0$ for any $j\ins\Z_{>0}$ using \cite[Thm.\,1]{bha}. Indeed, the latter theorem actually says that the roots of $b_f(s)$ supported at the origin are contained in $\Z\cap[n,2d{-}2]$ after multiplied by $-d$ (with $d:=\deg f$). Here $2d{-}2$ is closely related to $\Z_{\ges 2}$ in the assumption on $\al$. However, it does not seem easy to deduce Corollary~\hl{C3.2d}{3.2d} from Theorem~\hl{T2}{2} using \cite[Thm.\,1]{bha} (without using Corollary~\hl{C3.2b}{3.2b}) unless $D$ is tame. Indeed, it is not necessarily easy to show that $\DRl_X^{-1}\bl(\M^{\ssb}_{X,\log}(L),\naal\br)[n]$ is quasi-isomorphic to a $\D_X$-module. 
\par\htt{R3.2e}{}\msn
{\bf Remark~3.2e.} 
In the case $\al_k\eq\al$ ($\forall\,k$), Corollary~\hl{C3.2a}{3.2a} is compatible with Proposition~\hl{P1.5}{1.5} via the calculation of vanishing cycles using \cite[1.3]{BuSa}. Here $\om_p\wedge$ is given by $\al\,\ddd f\!/\!f\wedge$ (where $\al$ is not $\al{-}1$), and its action on $\A_{f,p}^{\ssb}$ vanishes by definition.
\par\htt{R3.2f}{}\msn
{\bf Remark~3.2f.} One can extend Corollary~\hl{C3.2d}{3.2d} to the case where the condition $\al_k\eq\al$ ($\forall\,k$) does not hold considering the action of $\D_X$ on $\mprod_k\,f_k^{\al_k}$ and modifying $\Tht_{f,\al-1}$ appropriately (where Hartogs theorem may be needed).
\bs\bs\htt{Ap}{}
\vbox{\centerline{\bf Appendix: Annihilator generated by vector fields}
\bsn
In this Appendix we give a proof of Theorem~\hl{TA}{A} below using an argument which is quite different from the proof of Theorem~\hl{T2}{2}, and is closely related to the theory of BS polynomials in \cite{Ka}, \cite{Bj}. In the case where the assumptions of Corollary~\hl{C2}{2} are satisfied, Theorem~\hl{TA}{A} together with this corollary implies a positive answer to a question in \cite[3.3]{To2} (using Corollary~\hl{C1.4}{1.4}), where $\al\eq{-}1$.}
\par\htt{TA}{}\msn
{\bf Theorem~A.} {\it Let $D$ be a reduced divisor on a complex manifold $X$ which is everywhere positively weighted homogeneous and tame. Then for $\al\ins\C$, the annihilator ${\rm Ann}_{\D_X}(f^{\al})$ of $f^{\al}$ in $\D_X$ is generated by logarithmic vector fields with $\OO$-linear terms if and only if $b_f(\al{-}j)\ne 0$ for any $j\ins\Z_{>0}$.}
\ms
(Here we use analytic $\D$-modules. The corresponding assertion for algebraic $\D$-modules follows by using the full faithfulness of $\OO_{X,x}$ over $\OO_{X_{\rm alg},x}$ for closed points $x\in X_{\rm alg}$.) We say that the annihilator ${\rm Ann}_{\D_X}(f^{\al})$ is {\it generated by logarithmic vector fields with $\OO$-linear terms\1} if it is generated over $\D_X$ by
$$\Tht_{f,\al}:=\bl\{\1\xi\mi\al\1\xi(f)/f\mid\xi\ins\Theta_X(-\log D)\br\}\q(\al\ins\C),$$
\par\nin where $\Theta_X(-\log D)$ denotes the sheaf of logarithmic vector fields \cite{SaK}, see also \hl{1.6}{1.6}. Note that $\Tht_{f,\al}$ is the specialization at $s\eq\al$ of
$$\Tht_f:=\bl\{\xi\mi s\1\xi(f)/f\mid\xi\ins\Theta_X(-\log D)\br\}\subset\D_X[s],$$
\par\nin which is contained in ${\rm Ann}_{\D_X[s]}(f^s)$. We have the isomorphisms as $\OO_X$-modules
\htt{A.1}{}
$$\Tht_f\eq\Tht_{f,\al}=\Theta_X(-\log D).
\leqno{\rm(A.1)}$$
\par\nin Theorem~\hl{TA}{A} is essentially a corollary of \cite[Thm.\,3.26]{Wa2} where $\Tht_f$ is used as is explained below.
\sk
For the proof of Theorem~\hl{TA}{A}, we recall some basic of BS polynomials, see \cite{Ka}, \cite{Bj}. For a holomorphic function $f$ on a complex manifold $X$, set $D:=f^{-1}(0)\sst X$, and
$$\aligned\Nc_f&:=\D_X[s]f^s\,\,(\subset\OO_X(*D)[s]f^s),\\ \Nc_{\al}&:=\Nc_f/(s{-}\al)\Nc_f\q(\al\ins\C).\endaligned$$
\par\nin There are canonical $\D_X$-linear surjective morphisms
\htt{A.2}{}
$$r_{\al}:\Nc_{\al}\onto\D_Xf^{\al}\,\,(\subset\,\OO_X(*D)f^{\al})\q(\al\ins\C).
\leqno{\rm(A.2)}$$
\par\nin We have also a $\D_X$-linear injective endomorphism $t:\Nc_f\into\Nc_f$ defined by
\htt{A.3}{}
$$t(P(s)f^s)=P(s{+}1)f^{s+1}\q(P(s)\ins\D_X[s]),
\leqno{\rm(A.3)}$$
\par\nin see \cite{Ka}. By definition $b_f(s)$ is the minimal polynomial of the action of $s$ on the holonomic $\D_X$-module $\Nc_f/t\Nc_f$ (which is induced by multiplication by $s$ on $\D_X[s]$). We assume $b_f(s)$ exists shrinking $X$ if necessary.
\sk
Since $ts\eq(s{+}1)t$, the morphism $t$ induces the $\D_X$-linear morphisms
\htt{A.4}{}
$$t_{\al}:\Nc_{\al}\to\Nc_{\al-1}\q(\al\ins\C),
\leqno{\rm(A.4)}$$
\par\nin (see also (\hl{A.6}{A.6}) below) together with the commutative diagram
\htt{A.5}{}
$$\begin{array}{cccccccc}\Nc_{\al}&\buildrel{\!\!t_{\al}}\over\longrightarrow&\Nc_{\al-1}&\buildrel{\!\!t_{\al-1}}\over\longrightarrow&\Nc_{\al-2}&\to&\cdots\\ \,\,\,\raise3mm\h{$\ontov$}\raise1.3pt\h{$\scriptstyle r_{\al}$}&\raise4mm\h{}\raise-2mm\h{}&\,\,\,\,\raise3mm\h{$\ontov$}\raise1.3pt\h{$\scriptstyle r_{\al-1}$}&&\,\,\,\,\raise3mm\h{$\ontov$}\raise1.3pt\h{$\scriptstyle r_{\al-2}$}\\ \D_Xf^{\al}&\buildrel{\iota_{\al}}\over\into&\D_Xf^{\al-1}&\buildrel{\iota_{\al-1}}\over\into&\D_Xf^{\al-2}&\into&\cdots\end{array}
\leqno{\rm(A.5)}$$
\par\nin Here the $\iota_{\al}$ are natural inclusions ($\al\ins\C$). It is easy to see that $\iota_{\al}$ is surjective if $b_f(\al{-}1)\nes0$, although the converse does not necessarily hold, see \cite[Ex.\,4.2]{rp}.
\sk
The following lemma is well known to specialists, see \cite[Lemma 6.21]{Ka2} and also a slightly weaker \cite[Prop.\,6.3.15]{Bj}. We note a short proof here for the convenience of the reader.
\par\htt{LA}{}\msn\vbox{\nin
{\bf Lemma~A.} {\it For $\al\ins\C$, the following conditions are equivalent\,$:$
\skn
{\rm(a)} $t_{\al}$ is injective,
\skn
{\rm(b)} $t_{\al}$ is surjective,
\skn
{\rm(c)} $b_f(\al{-}1)\nes0$.}}
\msn
{\it Proof.} This follows applying the snake lemma to the diagram below in view of the definition of $b_f(s)$ noted after (\hl{A.3}{A.3}).
\htt{A.6}{}
$$\begin{array}{cccccccc}&&&&0&\to&{\rm Ker}\,t_{\al}\\
&&&&\downarrow&&\downarrow\\
0&\to&\Nc_f&\buildrel{s-\al}\over\longrightarrow&\Nc_f&\to&\Nc_{\al}&\to0\\
&&\,\,\downarrow\!\raise1pt\h{$\scriptstyle t$}&&\,\,\downarrow\!\raise1pt\h{$\scriptstyle t$}&&\,\,\,\,\,\downarrow\!\raise1pt\h{$\scriptstyle t_{\al}$}\\
0&\to&\Nc_f&\buildrel{s-\al+1}\over\longrightarrow&\Nc_f&\to&\Nc_{\al-1}&\to0\\
&&\raise2mm\h{$\ontov$}&&\raise2mm\h{$\ontov$}&&\raise2mm\h{$\ontov$}\\
&&\Nc_f/t\Nc_f&\buildrel{s-\al+1}\over\longrightarrow&\Nc_f/t\Nc_f&\to&{\rm Coker}\,t_{\al}&\to0\end{array}
\leqno{\rm(A.6)}$$
\par\nin Indeed, $\Nc_f/t\Nc_f$ is a holonomic $\D_X$-module (having locally finite length) so that the action of $s{-}\al{+}1$ is injective if and only if it is surjective. This finishes the proof of Lemma~\hl{LA}{A}.
\sk
We have the following (which shows that the converse of \cite[Prop.\,6.2]{Ka} holds, see also \cite{Oa} for the case $f$ is a polynomial).
\par\htt{PA}{}\msn
{\bf Proposition~A.} {\it For $\al\ins\C$, the surjective morphism $r_{\al}:\Nc_{\al}\onto\D_Xf^{\al}$ in {\rm(\hl{A.2}{A.2})} is bijective if and only if $b_f(\al{-}j)\ne 0$ for any $j\ins\Z_{>0}$.}
\msn
{\it Proof.} It is enough to show that $b_f(\al{-}j)\nes0$ ($j\ins\Z_{>0})$ assuming the bijectivity of $r_{\al}$, since the converse is proved in \cite[Prop.\,6.2]{Ka}. The desired assertion, however, follows easily from Lemma~\hl{LA}{A} using the diagram (\hl{A.5}{A.5}) inductively. This finishes the proof of Proposition~\hl{PA}{A}.
\msn
{\it Proof of Theorem}~\hl{TA}{A}. The annihilator ${\rm Ann}_{\D_X[s]}f^s$ is generated by $\Tht_f$ over $\D_X[s]$, see \cite[Thm.\,3.26]{Wa2}. (Here the tameness assumption cannot be weakened as in Theorem~\hl{T2}{2}.) This implies that the annihilator of the image of $f^s$ in $\Nc_{\al}$ is generated by $\Tht_{f,\al}$ over $\D_X$. Indeed, we have the commutative diagram
\htt{A.7}{}
$$\begin{array}{cccccccc}
0&\to&\D_X[s]{\otimes}\Tht_f&\buildrel{\phi'}\over\to&\D_X[s]{\otimes}\Tht_f&\to&\D_X{\otimes}\Tht_{f,\al}&\to0\\ &\raise4.5mm\h{}&\downarrow&&\downarrow&&\,\,\,\downarrow\!\raise1pt\h{$\scriptstyle\rho$}\\
0&\to&\D_X[s]&\buildrel{\phi}\over\to&\D_X[s]&\to&\D_X&\to0\\
&&\raise2mm\h{$\ontov$}&&\raise2mm\h{$\ontov$}&&\raise2mm\h{$\ontov$}\\
0&\to&\,\Nc_f&\buildrel{\phi''}\over\to&\,\Nc_f&\to&\,{\rm Coker}\,\rho&\to0\end{array}
\leqno{\rm(A.7)}$$
\par\nin where $\phi',\phi,\phi''$ are induced by multiplication by $s\mi\al$ on $\D_X[s]$. We see that its rows and columns are exact using the snake lemma together with the isomorphisms in (\hl{A.1}{A.1}), where $\phi''$ is injective using the inclusion $\Nc_f\sst\OO_X(*D)[s]f^s$. From the diagram we thus get the canonical isomorphism
\htt{A.8}{}
$$\Nc_{f,\al}=\D_X/\D_X\Tht_{f,\al}.
\leqno{\rm(A.8)}$$
\par\nin So the assertion follows from Proposition~\hl{PA}{A}. This completes the proof of Theorem~\hl{TA}{A}.

\end{document}